\documentclass{article}
\usepackage{eepic,amsthm,amssymb,amsmath}
\usepackage{amssymb}
\usepackage{epsfig,array,picinpar,psfrag,graphics}
\usepackage[mathscr]{eucal}

 \allinethickness{0.4mm}

\newtheorem{prop}{Proposition}

\newtheorem{theorem}{Theorem}
\theoremstyle{definition}
\newtheorem{algorithm}{Algorithm}
\newtheorem{definition}{Definition}

\newcommand{\Aut}{\mathop{\rm Aut}\nolimits}

\newcommand{\Stg}{\mathop{\mathsf{St}}\nolimits_G}

\newcommand\Z{{\mathbb Z}}

\newcommand\C{{\mathbb C}}
\renewcommand\H{{\mathcal H}}

\newcommand\Sym{{\mathsf{Sym}}}
\newcommand\Isom{{\mathsf{Isom}}}

\newcommand\tree{{\mathcal{T}}}
\newcommand\btree{{\partial\tree}}
\newcommand\St{{\mathsf{St}}}

\newcommand\G{{\mathcal{G}}}
\newcommand\B{{\mathcal{B}}}

\newcommand{\lims}[1][G]{\mathscr{J}_{#1}}
\newcommand{\gSr}{\Gamma}

\title{Groups generated by $3$-state automata \\
over a $2$-letter alphabet, I}

\author{Ievgen Bondarenko, Rostislav Grigorchuk, Rostyslav Kravchenko,\\ Yevgen Muntyan, Volodymyr
Nekrashevych,\\ Dmytro Savchuk and Zoran \v{S}uni\'{c}\thanks{All
authors were partially supported by at least one of the NSF grants
DMS-308985, DMS-456185, DMS-600975, and DMS-605019} }

\begin{document}
\maketitle

\begin{abstract}
An approach to a classification of groups generated by 3-state
automata over a 2-letter alphabet and the current progress in this
direction are presented. Several results related to the whole class
are formulated. In particular, all finite, abelian, and free groups
are classified. In addition, we provide detailed information and
complete proofs for several groups from the class, with the
intention of showing the main methods and techniques used in the
classification.
\end{abstract}

\section{Introduction}
Groups generated by finite automata (groups of automata or automaton
groups) were formally introduced at the beginning of
1960's~\cite{horejs:automata}, but more substantial work on this
remarkable class of groups started only in 1970's after
Aleshin~\cite{aleshin:burnside} confirmed a conjecture by
Glushkov~\cite{glushkov:ata} that these groups could be used to
study problems of Burnside type (note that groups of automata should
not be confused with automatic groups as described
in~\cite{epstein-al:wp}). It was observed in 1960's and 1970's that
groups of automata are closely related to iterated wreath products
(pioneering work in this direction is due to
Kaloujnin~\cite{kalouj:psubgr1}) and that the theory of such groups
could be studied by using the language of tables developed by
Kaloujnin~\cite{kalouj:psubgr2} and
Sushchanskii~\cite{sushchanskii:burnside}.

Even more intensive study of groups of finite automata started in
the beginning of 1980's after the development of some new ideas such
as self-similarity, contracting properties, and a geometric
realization as groups acting on rooted trees. These developments
allowed for elegant constructions of Burnside
groups~\cite{grigorchuk:burnside,gupta-s:burnside,gupta_s:pgroups}
and pushed the study of groups of automata in many directions:
analysis~\cite{grigorchuk:gdegree,erschler:subexponential},
geometry~\cite{bartholdi-g-n:fractal},
probability~\cite{bartholdi-v:basilica,erschler:subexponential,abert-v:dimension},
dynamics~\cite{bartholdi_g:spectrum,grigorchuk-z:l2}, formal
languages~\cite{holt-r:word}, etc.

\begin{figure}
\begin{center}
\epsfig{file=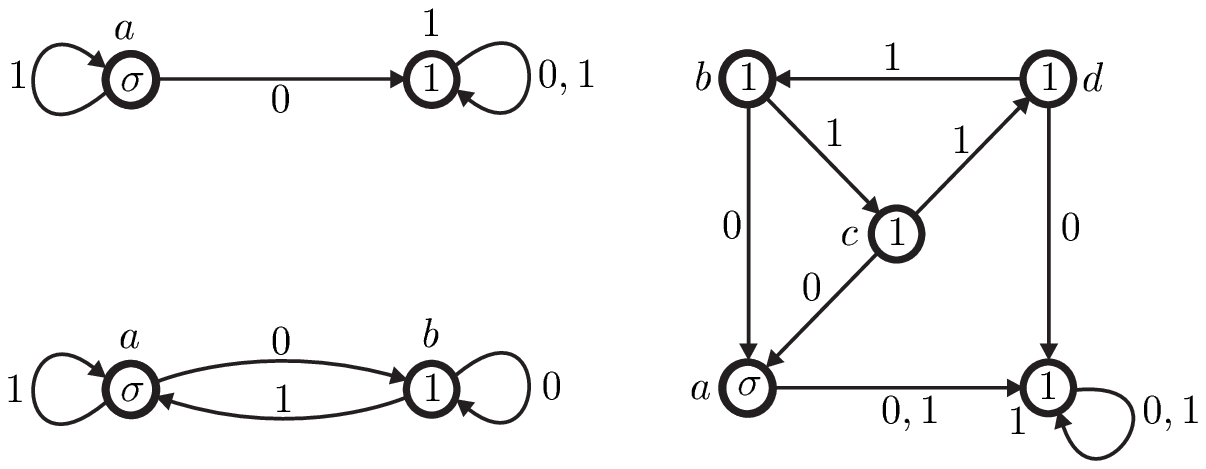,height=110pt}
\end{center}
\caption{Adding machine, lamplighter automaton and automaton
generating $\G$\label{3automata}}
\end{figure}

Two well known and important problems were solved using groups of
automata in the early 1980's, namely Milnor
Problem~\cite{milnor:problem} on intermediate growth and Day
Problem~\cite{day:amenable} on amenability. A 5-state automaton
constructed in~\cite{grigorchuk:burnside} (on the right in
Figure~\ref{3automata}) generates a 2-group, denoted $\G$. It was
shown in~\cite{grigorchuk:gdegree} that $\G$ has intermediate growth
(between polynomial and exponential). This led to construction of
other examples of this
type~\cite{fabrykovski-g:growth2,bartholdi-s:wpg} and also made
important contribution to and impact on the theory of invariant
means on
groups~\cite{greenleaf:b-means,wagon:b-paradox,paterson:b-amenability}
initiated by von Neumann~\cite{vneumann:masses} by providing an
example of amenable, but not elementary amenable group (in the sense
of Day~\cite{day:amenable}).

Among the most interesting newer developments is the spectral theory
of groups generated by finite automata and graphs associated to such
groups~\cite{bartholdi_g:spectrum,grigorchuk-z:l2}. For instance,
automaton groups provided first examples of regular graphs realized
as Schreier graphs of groups for which the spectrum of the
combinatorial Laplacian is a Cantor set~\cite{bartholdi_g:spectrum}.
Further, the realization of the lamplighter group $\Z\wr C_2$ as
automaton group (bottom left in Figure~\ref{3automata}) was crucial
in the proof that this group has a pure point spectrum (with respect
to a system of generators related to the states of the automaton)
and thus has discrete spectral measure, which was completely
described~\cite{grigorchuk-z:l2}. This, in turn, led to a
construction~\cite{grigorchuk-al:atiyah} of a 7-dimensional closed
manifold with non-integer third $L^2$-Betti number providing a
counterexample to the Strong Atiyah Conjecture~\cite{luck:b-l2}.

\begin{figure}
\begin{center}
\epsfig{file=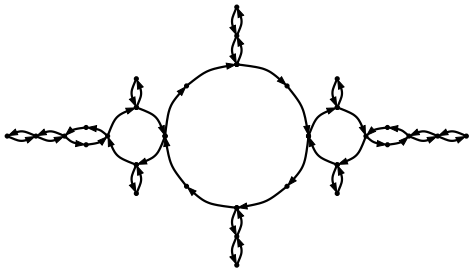,height=100pt}
\end{center}
\caption{Schreier graph of the action of Basilica group $\B$ on the
5-th level of the tree\label{b-schreier5}}
\end{figure}

Another fundamental recent discovery is the relation of groups of
automata to holomorphic
dynamics~\cite{bartholdi-g-n:fractal,nekrash:self-similar}. Namely,
it is shown that to every rational map $f: \hat\C \to \hat\C$ on the
Riemann Sphere with finite postcritical set one can associate a
finite automaton generating a group, denoted $IMG(f)$ and called
iterated monodromy group of $f$. The geometry and the topology of
the Schreier graphs of $IMG(f)$ is closely related to the geometry
of the Julia set of $f$. Figure~\ref{b-schreier5} depicts a Schreier
graph associated to an automaton group, denoted $\B$ and called
Basilica group. Its reminiscence to the Julia set of the map $z
\mapsto z^2 -1$ is related to the fact that $\B$ is the iterated
monodromy group $IMG(z^2-1)$ of the holomorphic map $z \mapsto z^2
-1$~\cite{bartholdi-g-n:fractal,nekrash:self-similar}). Groups of
automata represent the basis of the theory of self-similar groups
and actions~\cite{nekrash:self-similar} and are related to the study
of Belyi polynomials and dessins d'enfants of
Grothendieck~\cite{pilgrim:dessins}. The use of iterated monodromy
groups was crucial in the recent solution of Hubbard's Twisted
Rabbit Problem in~\cite{bartholdi_n:rabbit}.

One of the most important developments in the theory of automaton
groups is the introduction of groups with branch structure,
providing a link to just-infinite
groups~\cite{grigorchuk:jibg,wilson:jipg,bartholdi-g-s:branch} and
groups of finite width~\cite{bartholdi-g:lie}. In particular, a
problem suggested by Zelmanov was solved by using the profinite
completion of $\G$~\cite{bartholdi-g:lie}. The problem of Gromov on
uniformly exponential growth was solved recently by using branch
automaton groups~\cite{wilson:nonuniform}. An unexpected link of
groups of automata and their profinite completions to Galois theory
was found by R.~Pink (private communication) and Aitken, Hajir and
Maire~\cite{aitken-h-m:iterated}, while N.~Boston
\cite{boston:survey,boston:reducing} related branch groups of
automata to Fontaine-Mazur Conjecture and other problems in number
theory. The class of branch groups is also a new source for
infinitely presented groups, for which the presentation can be
written in a recursive form
(see~\cite{lysionok:presentation,sidki:presentation,grigorch_ss:img}).

A recent observation~\cite{grigorchuk-s:hanoi-crm} is that automaton
groups and their Schreier graphs stand behind the famous Hanoi
Towers Problem (see~\cite{hinz:towers}) and some of its
generalizations~\cite{stockmeyer:variations}.

There are indications that spectral properties of groups generated
by finite automata could be used in the study of Kaplansky
Conjecture on Idempotents (and thus also Baum-Connes
Conjecture~\cite{connes:b-ncg} and Novikov
Conjecture~\cite{ferry-r-r:novikovconj}), Dixmier Unitarizability
Problem~\cite{dixmier:problem,pisier:are}, and for construction of
new families of expanders, and perhaps even Ramanujan,
graphs~\cite{lubotzky:expander-book}.


In this article we are going to describe some progress which was
achieved during the last few years in the problem of classification
of automaton groups.

Two important characteristics of an automaton are the cardinalities
$m$ and $n$ of the set of states and the alphabet, respectively, and
the pair $(m,n)$ is a natural measure of complexity of an automaton
and of the group it generates.

The groups of complexity $(2,2)$ are classified
\cite{gns00:automata} and there are only $6$ such groups (see
Theorem~\ref{thm_class2states} in Section~\ref{sec_proofs} here).
The problem of classification of $(3,2)$ groups or $(2,3)$ groups is
much harder.

The current text represents the progress being made by the research
group at Texas A\&M University over the last few years toward
classification of $(3,2)$ groups.

The total number of invertible automata of complexity $(3,2)$ is
$2^3 \cdot 3^6 = 5832$. However, the number of non-isomorphic groups
generated by these automata is much smaller.

\begin{theorem}
There are no more than $124$ pairwise non-isomorphic groups of
complexity $(3,2)$.
\end{theorem}

The proof of this theorem is too long to be presented here (even the
list of all groups takes a lot of space).

Instead, we have chosen for this article a set of $24$ automata
generating $20$ groups (among the most interesting in this class, in
our opinion), which we list in the form of a table. The table
provided here is a part of the table listing the whole set of $124$
groups. We keep the numeration system from the whole table (the rule
for numeration is explained in Section~\ref{sec_approach}).

Major results obtained for the whole family are the following
theorems. The numbers in the brackets indicate the numbers of
corresponding automata in the class.

\begin{theorem}
There are $6$ finite groups in the class: $\{1\}$ [1], $C_2$ [1090],
$C_2\times C_2$ [730], $D_4$ [847], $C_2\times C_2\times C_2$ [802]
and $D_4\times C_2$ [748].
\end{theorem}

\begin{theorem}
There are $6$ abelian groups in the class: $\{1\}$ [1], $C_2$
[1090], $C_2\times C_2$ [730], $C_2\times C_2\times C_2$ [802],
$\mathbb Z$ [731] and $\mathbb Z^2$ [771].
\end{theorem}

Note that there are also virtually abelian groups in this class
(having $\Z$, $\Z^2$ [2212], $\Z^3$ [752] or $\Z^5$ [968] as
subgroups of finite index).

\begin{theorem}
\label{thm:class_free} The only free non-abelian group in the class
is the free group of rank 3 generated by the Aleshin-Vorobets
automaton [2240]. Moreover, the isomorphism class of this automaton
group coincides with its equivalence class under symmetry.
\end{theorem}

The definition of symmetric automata is given in
Section~\ref{sec_approach}.

\begin{theorem}
There are no infinite torsion groups in the class.
\end{theorem}

We do not provide the complete proofs of these theorems (by the
reason explained above). Instead, we give here some information
about each of the chosen groups and include proofs of most facts.

Properties that are in our focus are the contracting property,
self-replication, torsion, relations (we list the relators up to
length $10$), rank of quotients of the stabilizers series, shape of
the related Schreier graphs. The article is organized as follows. We
start with a general information about rooted trees and their
automorphisms. Then we provide quick introduction to the theory of
automaton groups. We continue with the definition of Schreier graphs
and explain how they naturally appear for the actions on rooted
trees. Then we list $24$ automata generating $20$ groups together
with some of their properties. In the last section we give proofs of
many facts related to the groups in the list.

The last part also contains some more general results (such as an
algorithm detecting transitivity of an element and a criterion for
group transitivity on the binary tree).

We recommend the articles~\cite{gns00:automata} and the book
\cite{nekrash:self-similar} to the reader who is interested in
becoming more familiar with automaton groups.


\section{Regular rooted tree automorphisms and self-similarity}

Let $d \geq 2$ be fixed and let $X$ be the alphabet
$X=\{0,1,\dots,d-1\}$. The set of words $X^*$ over $X$ (the free
monoid over $X$) can be given the structure of a \emph{regular
rooted labeled $d$-ary tree} $\tree$ in which the empty word
$\emptyset$ is the \emph{root}, the \emph{level} $n$ in $\tree$
consists of the words of length $n$ over $X$ and every vertex $v$
has $d$ children, labeled by $vx$, for $x \in X$. Denote by
$\Aut(\tree)$ the group of automorphisms of $\tree$. Let $f$ be an
automorphism in $\Aut(\tree)$. Any such automorphism can be
decomposed as
\begin{equation}\label{decomposition}
  f = \alpha_f(f_0,\dots,f_{d-1})
\end{equation}
where $f_x$, for $x \in X$, are automorphisms of $\tree$ and
$\alpha_f$ is a permutation of $X$. The automorphisms $f_x$ (also
denoted by $f|_x$), $x \in X$, are called (the first level)
\emph{sections} of $f$ and each one acts as an automorphism on the
subtree $\tree_x$ hanging below the vertex $x$ in $\tree$ consisting
of the words in $X^*$ that start with $x$ (any such subtree is
canonically isomorphic to the whole tree). The action of $f$ is
decomposed in two steps. First the $d$-tuple $(f_0,\dots,f_{d-1})$
acts on the $d$ subtrees hanging below the root, and then the
permutation $\alpha_f$, called the \emph{root permutation} of $f$,
permutes these $d$ subtrees. Thus the action of $f$ from
\eqref{decomposition} on $X^*$ is given by
\begin{equation} \label{fxw}
 f(xw) = \alpha_f(x) f_x(w),
\end{equation}
for $x$ a letter in $X$ and $w$ a word over $X$. Further iterations
of the decomposition (\ref{decomposition}) yield the second level
sections $f_{xy}=(f_x)_y$, $x,y \in X$, and so on. Algebraically, we
have
\begin{equation}
\label{eqn:iter_wreath}
\begin{array}{r}
\Aut(\tree) = \Sym(X) \ltimes(\Aut(\tree), \dots, \Aut(\tree))=\\
\Sym(X)\ltimes\Aut(\tree)^X   = \Sym(X)\wr\Aut(\tree),
\end{array}
\end{equation}
where $\wr$ is the \emph{permutational wreath product} in which the
coordinates of $\Aut(\tree)^X$ are permuted by $\Sym(X)$.

Iterations of the decomposition~\eqref{eqn:iter_wreath} show that
$\Aut(\tree)$ has the structure of an iterated wreath product
$\Aut(\tree) = \Sym(X) \wr(\Sym(X)\wr(\Sym(X)\wr\dots))$. Thus
$\Aut(\tree)$ is a pro-finite group and in particular, all of its
subgroups are residually finite. An obvious and natural sequence of
normal subgroups of finite index intersecting trivially is the
sequence of level stabilizers. The $n$-th \emph{level stabilizer}
$\St_{G}(n)$ of a group $G \leq \Aut(\tree)$ consists of those tree
automorphisms in $G$ that fix the vertices in $\tree$ up to level
$n$. The group $\Aut(\tree)$ is obviously an uncountable object. We
are interested in finitely generated subgroups of $\Aut(\tree)$ that
exhibit some important features of $\Aut(\tree)$. One such feature
is self-similarity.

\begin{definition}
A group $G$ of tree automorphisms is \emph{self-similar} if, for
every $g$ in $G$ and a letter $x$ in $X$ there exists a letter $y$
in $X$ and an element $h$ in $G$ such that
\[ g(xw) = yh(w), \]
for all words $w$ over $X$.
\end{definition}

Another way to express self-similarity of a group $G$ of tree
automorphisms is to say that every section $g_x$ of every element
$g$ in $G$ is again an element of $G$. The full tree automorphism
group $\Aut(\tree)$ is clearly self-similar (see (\ref{fxw})). A
self-similar group $G$ embeds in the permutational wreath product
$\Sym(X)\wr G = \Sym(X)\ltimes G^X$ by
\begin{equation}
\label{eqn:wr_rec} g \mapsto \alpha_g(g_0,g_1,\dots,g_{d-1}).
\end{equation}

\section{Definition of automaton groups}\label{definition}

Consider a finite system of recursive relations
\begin{equation} \label{f^0}
\left\{
 \begin{array}{ccc}
 f^{(1)} &= & \alpha_1 \left(f^{(1)}_0,f^{(1)}_1,\dots,f^{(1)}_{d-1}\right),  \\
 \dots  \\
 f^{(k)} &= & \alpha_k \left(f^{(k)}_0,f^{(k)}_1,\dots,f^{(k)}_{d-1}\right),
 \end{array}
\right.
\end{equation}
where each symbol $f^{(i)}_j$, $i=1,\dots,k$, $j=0,\dots,d-1$, is
equal to one of the symbols $f^{(1)}, \dots, f^{(k)}$ and
$\alpha_1,\dots,\alpha_k \in \Sym(X)$. The system (\ref{f^0}) has a
unique solution in $\Aut(\tree)$. The action of $f^{(i)}$ on $\tree$
is given recursively by $ f^{(i)}(xw) = \alpha_i(x) f^{(i)}_x (w)$.
The group generated by the automorphisms $f^{(1)}, \dots, f^{(k)}$
is finitely generated self-similar group of automorphisms of
$\tree$. This group can be described by a finite invertible
automaton (just called automaton in the rest of the article).

\begin{definition} \label{automaton}
A \emph{finite invertible automaton} $A$ is a $4$-tuple
$A=(Q,X,\rho,\tau)$ where $Q$ is a finite set of \emph{states}, $X$
is a finite \emph{alphabet} of cardinality $d \geq 2$, $\rho:Q
\times X \to X$ is a map, called \emph{output map}, $\tau:Q \times X
\to Q$ is a map, called \emph{transition map}, and for each state
$q$ in $Q$, the restriction $\rho_q: X \to X$ given by
$\rho_q(x)=\rho(q,x)$ is a permutation, i.e. $\rho_q \in \Sym(X)$.
\end{definition}

The automaton $A=(Q,X,\rho,\tau)$ reads words from $X^*$ and
provides output words that are also in $X^*$. The behavior is
encoded in the output and transition maps. An \emph{initial
automaton} $A_q$ is just an automaton $A$ with a distinguished state
$q \in Q$ selected as an initial state. We first informally describe
the action of the initial automaton $A_q$ on $X^*$. The automaton
starts at the state $q$, reads the first input letter $x_1$, outputs
the letter $\rho_q(x_1)$ and changes its state to $q_1=\tau(q,x_1)$.
The rest of the input word is handled by the new state $q_1$ in the
same fashion (in fact it is handled by the initial automaton
$A_{q_1}$). Formally, the action of the states of the automaton $A$
on $X^*$ can be described by extending the output function $\rho$ to
a function $\rho:Q \times X^* \to X^*$ recursively by
\begin{equation} \label{rhoext}
\rho(q,xw) = \rho(q,x) \rho(\tau(q,x),w)
\end{equation}
for all states $q$ in $Q$, letters $x \in X$ and words $w$ over $X$.
Then the action of the initial automaton $A_q$ is defined by
$A_q(u)=\rho(q,u)$, for words $u$ over $X$. In fact, (\ref{rhoext})
shows that each initial automaton $A_q$, $q \in Q$, defines a tree
automorphism, denoted by $q$, defined by
\begin{equation} \label{qxw}
 q(xw) = \alpha_q(x) q_x(w),
\end{equation}
where the section $q_x$ is the state $\tau(q,x)$ and the root
permutation $\alpha_q$ is the permutation $\rho_q$.

\begin{definition}
Given an automaton $A=(Q,X,\rho,\tau)$, the group of tree
automorphisms generated by the states of $A$ is denoted by $G(A)$
and called the \emph{automaton group} defined by $A$. The generating
set $Q$ is called the \emph{standard} generating set of $G(A)$.
\end{definition}

\emph{Boundary} of the tree $\tree$, denoted $\btree$, is the set
$X^\omega$ of words over $X$ that are infinite to the right
(infinite geodesic rays in $\tree$ starting at the root). It has a
natural metric (infinite words are close if they agree on long
finite prefixes) and the group of isometries $\Isom(\btree)$ is
canonically isomorphic to $\Aut(\tree)$. Thus the action of the
automaton group $G(A)$ on $\tree$ can be extended to an isometric
action on $\btree$. In fact, \eqref{rhoext} and \eqref{qxw} are
valid for infinite words $w$ as well.

An automaton $A$ can be represented by a labeled directed graph,
called Moore diagram, in which the vertices are the states of the
automaton, each state $q$ is labeled by its own root permutation
$\alpha_q$ and, for each pair $(q,x) \in Q \times X$, there is an
edge from $q$ to $q_x=\tau(q,x)$ labeled by $x$. For example, the
5-state automaton in the right half of Figure~\ref{3automata}
generates the group $\G$ mentioned in the introduction ($\sigma$
denotes the permutation exchanging $0$ and $1$). The two 2-state
automata given on the left of Figure~\ref{3automata} are the so
called \emph{adding machine} (top), which generates the infinite
cyclic group $\Z$ and the \emph{lamplighter automaton} (bottom)
generating $L_2=\Z\wr C_2$. Recursion relations of type (\ref{f^0})
for the adding machine and the lamplighter automaton are given by
\begin{alignat}{4}
  &a = \sigma &&(1,a)  \qquad\qquad\qquad\qquad\qquad\qquad &&a = \sigma &&(b,a) \notag\\
  &1 =        &&(1,1)                                       &&b = &&(b,a), \label{lampl}
\end{alignat}
respectively.

Various classes of automaton groups deserve special attention. An
automaton group $G=G(A)$ is \emph{contracting} if there exist
constants $\kappa$, $C$, and $N$, with $0 \leq \kappa < 1$, such
that $|g_v| \leq \kappa |g| + C$, for all vertices $v$ of length at
least $N$ and $g \in G$ (the length is measured with respect to the
standard generating set $Q$). For sufficiently long elements $g$
this means that the length of its sections at vertices on levels
deeper than $N$ is strictly shorter than the length of $g$. This
length shortening leads to an equivalent definition of a contracting
group. Namely, a group $G$ of tree automorphisms is contracting if
there exists a finite set $\mathcal N\subset G$, such that for every
$g\in G$, there exists $N>0$, such that $g_v\in\mathcal N$ for all
vertices $v\in X^*$ of length not shorter than $N$. The minimal set
$\mathcal N$ with this property is called the \emph{nucleus} of $G$.
The contraction property is a key feature of various inductive
arguments and algorithms involving the decomposition
$g=\alpha_g(g_0,\dots,g_{d-1})$.

Another important class is the class of automaton groups of
\emph{branch type}. Branch groups arise as one of the
three~\cite{grigorchuk:jibg} possible types of just-infinite groups
(infinite groups for which all proper homomorphic images are
finite). Every infinite, finitely generated group has a
just-infinite image. Thus if a class of groups $\mathcal{C}$ is
closed under homomorphic images and if it contains infinite,
finitely generated examples then it contains just-infinite examples.
Such examples are, in a sense, minimal infinite examples in
$\mathcal{C}$. For example, $\G$ is a branch automaton group that is
a just-infinite 2-group. i.e., it is an infinite, finitely
generated, torsion group that has no proper infinite quotients.
Also, the Hanoi Towers group~\cite{grigorchuk-s:hanoi-crm} and the
iterated monodromy group $IMG(z^2+i)$~\cite{grigorch_ss:img} are
branch groups, while $\B=IMG(z^2-1)$ is not a branch group, but only
weakly branch (for definitions
see~\cite{grigorchuk:jibg,bartholdi-g-s:branch}).

The class of \emph{polynomially growing automata} was introduced by
Sidki in~\cite{sidki:pol}, where it is proved that no group $G(A)$
defined by such an automaton contains free subgroups of rank 2.
Moreover, for a subclass of so called \emph{bounded automata} it is
known that the corresponding groups are
amenable~\cite{bartholdi-k-n-v:bounded} (this class of automata, for
instance, includes the automata generating $\G$, $\B$ and Hanoi
Towers group on $3$ pegs, but not for more pegs).

Finally, self-replicating groups play an important role. A
self-similar group $G$ is called \emph{self-replicating} if, for
every vertex $u$, the homomorphism $\varphi_u:\Stg(u) \to G$ from
the stabilizer $u$ in $G$ to $G$, given by $\varphi(g)=g_u$, is
surjective. This condition is usually easy to check and, together
with transitivity of the action on level 1, it implies transitivity
of the action on all levels. Another way to show that a group of
automorphisms of the binary tree is level transitive is to use
Proposition~\ref{prop:transitivity}.

\section{Limit spaces, Schreier graphs and iterated monodromy groups}

Let us fix some self-similar contracting group acting on $X^{*}$ by
automorphisms. Denote by $X^{-\omega}$ the space of left infinite
sequences over $X$.

\begin{definition}
Two elements $\ldots x_3x_2x_1, \ldots y_3y_2y_1\in X^{-\omega}$ are
said to be \emph{asymptotically equivalent} with respect to the
action of the group $G$, if there exist a finite set $K\subset G$
and a sequence $\{g_k\}_{k=1}^\infty$ of elements in $K$ such that
\[
g_k(x_kx_{k-1}\ldots x_2x_1)=y_ky_{k-1}\ldots y_2y_1
\]
for every $k\geq 1$.
\end{definition}

The asymptotic equivalence is an equivalence relation. Moreover,
sequences $\ldots x_2x_1, \ldots y_2y_1\in X^{-\omega}$ are
asymptotically equivalent if and only if there exists a sequence
$\{h_k\}$ of the elements in the nucleus of $G$ such that
$h_k(x_k)=y_k$ and $h_k|_{x_k}=h_{k-1}$, for all $k\geq 1$.

\begin{definition}
The quotient space $\lims$ of the topological space $X^{-\omega}$ by
the asymptotic equivalence relation is called the \emph{limit space}
of the self-similar action of $G$.
\end{definition}

The limit space $\lims$ is metrizable and finite-dimensional. If the
group $G$ is finitely-generated and level-transitive, then the limit
space $\lims$ is connected.

The last decade witnessed a shift in the attention payed to the
study of Schreier graphs. Let $G$ be a group generated by a finite
set $S$ and let $G$ act on a set $Y$. The \emph{Schreier graph} of
the action $(G,Y)$ is the graph $\gSr(G,S,Y)$ with set of vertices
$Y$ and set of edges $S\times Y$, where the arrow $(s,y)$ starts in
$y$ and ends in $s(y)$. If $y\in Y$ then the Schreier graph
$\gSr(G,S,y)$ of the action of $G$ on the $G$-orbit of $y$ is called
\emph{orbital Schreier graph}.

Let $G$ be a subgroup of $\Aut(\tree)$ generated by a finite set $S$
(not necessary self-similar). The levels $X^n$, $n \geq 0$, are
invariant under the action of $G$ and we can consider the Schreier
graphs $\gSr_n(G,S)=\gSr(G,S,X^n)$. Let $\omega=x_1x_2x_3\ldots\in
X^{\omega}$. Then the pointed Schreier graphs
$(\gSr_n(G,S),x_1x_2\ldots x_n)$ converge in the local topology
(topology defined in~\cite{grigorchuk:gdegree}) to the pointed
orbital Schreier graph $(\gSr(G,S,\omega),\omega)$.

The limit space of a finitely generated contracting self-similar
group $G$ can be viewed as a hyperbolic boundary in the following
way. For any given finite generating system $S$ of $G$ define the
self-similarity graph $\Sigma(G,S)$ as the graph with set of
vertices $X^{*}$ in which two vertices $v_1,v_2\in X^{*}$ are
connected by an edge if and only if either $v_i=xv_j$, for some
$x\in X$ (vertical edges), or $s(v_i)=v_j$ for some $s\in S$
(horizontal edges). If the group is contracting then the
self-similarity graph $\Sigma(G,S)$ is Gromov-hyperbolic and its
hyperbolic boundary is homeomorphic to the limit space $\lims$. The
set of horizontal edges of $\Sigma(G,S)$ spans the disjoint union of
all Schreier graphs $\gSr_n(G,S)$.Thus, the Schreier graphs
$\gSr_n(G,S)$ in some sense approximate the limit space $\lims$ of
the group $G$. Moreover, for many examples of self-similar
contracting groups there exists a sequence of numbers $\lambda_n$
such that the metric spaces $(\gSr_n,d(\cdot,\cdot)/\lambda_n)$,
where $d$ is the combinatorial metric on the graph, converge in the
Gromov-Hausdorff metric to the limit space of the group.

We recall the definition and basic properties of iterated monodromy
groups (IMG). Let $\mathcal{M}$ be a path connected and locally path
connected topological space and let $\mathcal{M}_1$ be its open path
connected subset. Let $f:\mathcal{M}_1\rightarrow\mathcal{M}$ be a
$d$-fold covering. By $f^n$ we denote the $n$-th iteration of the
map $f$. The map $f^n:\mathcal{M}_n\rightarrow\mathcal{M}$, where
$\mathcal{M}_n = f^{-n}(\mathcal{M})$, is a $d^n$-fold covering.

Choose an arbitrary base point $t\in\mathcal{M}$. Let $\tree_t$ be
the disjoint union of the sets $f^{-n}(t), n\geq 0$ (these sets are
note necessarily disjoint by themselves). The set of pre-images
$\tree_t$ has a natural structure of a rooted $d$-ary tree with root
$t\in f^{-0}(t)$ in which every vertex $z\in f^{-n}(t)$ is connected
to the vertex $f(z)\in f^{-n+1}(t)$, $n \geq 1$. The fundamental
group $\pi_1(\mathcal{M},t)$ acts naturally on every set $f^{-n}(t)$
and, in fact, acts by automorphisms on $\tree_t$.

\begin{definition}
\emph{Iterated monodromy group} $IMG(f)$ of the covering $f$ is the
quotient of the fundamental group $\pi_1(\mathcal{M},t)$ by the
kernel of its action on $\tree_t$.
\end{definition}

It is proved in~\cite{nekrash:self-similar} that all iterated
monodromy groups are self-similar. This fact provides a connection
between holomorphic dynamics and groups generated by automata.

\begin{theorem}
The iterated monodromy group of a sub-hyperbolic rational function
is contracting and its limit space is homeomorphic to the Julia set
of the rational function.
\end{theorem}

In particular, the sequence of Schreier graphs $\gSr_n$ of the
iterated monodromy group of a sub-hyperbolic rational function can
be drawn on the Riemann sphere in such a way that they converge in
the Hausdorff metric to the Julia set of the function.

Schreier graphs also play a role in computing the spectrum of the
Markov operator $M$ on the group. Namely, given a group $G$
generated by a finite set $S=\{s_1,s_2,\dots,s_k\}$, acting on a
tree $X^*$ there is a natural unitary representation of $G$ in the
space of bounded linear operators $\H=B(L_2(X^\omega))$ given by
$\pi_g(f)(x)=f(g^{-1}x)$.

The Markov operator
$M=\frac1{2k}(\pi_{s_1}+\dots+\pi_{s_k}+\pi_{s_1^{-1}}+\dots+\pi_{s_k^{-1}})$
corresponding to this unitary representation plays an important
role. The spectrum of $M$ for a self-similar group $G$ is
approximated by the spectra of finite dimensional operators arising
from the action of $G$ on the levels of the tree $X^*$. For more on
this see~\cite{bartholdi_g:spectrum}.

Let $\H_n$ be a subspace of $\H$ spanned by the $|X|^n$
characteristic functions $f_v, v\in X^n$, of the cylindrical sets
corresponding to the $|X|^n$ vertices on level $n$. Then $\H_n$ is
invariant under the action of $G$ and $\H_n\subset \H_{n+1}$. Denote
by $\pi^{(n)}_g$ the restriction of $\pi_g$ on $\H_n$. Then
\[M_n=\frac1{2k}(\pi^{(n)}_{s_1}+\pi^{(n)}_{s_2}+\dots+\pi^{(n)}_{s_k}+\pi^{(n)}_{s_1^{-1}}+\pi^{(n)}_{s_2^{-1}}+\dots+\pi^{(n)}_{s_k^{-1}})\]
are finite dimensional operators, whose spectra converge to the
spectrum of $M$ in the sense
\[sp(M)=\overline{\bigcup_{n\geq0}sp(M_n)}.\]

If $P$ is the stabilizer of an infinite word from $X^\omega$, then
one can consider the Markov operator $M_{G/P}$ on the Schreier graph
of $G$ with respect to $P$. The following fact is observed
in~\cite{bartholdi_g:spectrum} and can be applied to compute the
spectrum of Markov operator on the Cayley graph of a group in case
if $P$ is small.

\begin{theorem}
If $G$ is amenable or the Schreier graph $G/P$ (the Schreier graph
of the action of $G$ on the cosets of $P$) is amenable then
$sp(M_{G/P})=sp(M)$.
\end{theorem}


\section{Approach to a classification of groups generated by 3-state
automata over a 2-letter alphabet} \label{sec_approach}

The next three sections are devoted to the groups generated by
3-state automata over the 2-letter alphabet $X=\{0,1\}$. Fix
$\{\textbf1,\textbf2,\textbf3\}$ as the set of states. Every $(3,2)$
automaton is given by
\[
 \left\{
 \begin{array}{l}
  \textbf1=\sigma^{a_{11}}(a_{12},a_{13}), \\
  \textbf2=\sigma^{a_{21}}(a_{22},a_{23}),  \\
  \textbf3=\sigma^{a_{31}}(a_{32},a_{33}),
 \end{array}
 \right.
\]
where $a_{ij}\in \{\textbf1,\textbf2,\textbf3\}$, for $j\neq1$,
$a_{i1}\in \{0,1\}$, $i=1,2,3$, and $\sigma=(01)\in\Sym(X)$. A
number is assigned to the automaton above by the following formula
\[
\begin{array}{l}
\mathop{\rm Number}\nolimits(\mathcal A)= \\
\qquad\qquad(a_{12}-1)+3(a_{13}-1)+9(a_{22}-1)+27(a_{23}-1)+81(a_{32}-1)+\\
 \qquad\qquad 243(a_{33}-1)+729(a_{11}+2a_{21}+4a_{31})+1
\end{array}.
\]

Thus every $(3,2)$ automaton obtains a unique number in the range
from $1$ to $5832$. The numbering of the automata is induced by the
lexicographic ordering of all automata in the class. The automata
numbered $1$ through $729$ act trivially on the tree and generate
the trivial group. The automata numbered $5104$ through $5832$
generate the group $C_2$ of order $2$, because every element in any
of these groups is either trivial, or changes all letters in any
word over $X$. Therefore the ``interesting'' automata have numbers
$730$ through $5103$.

Denote by $\mathcal A_n$ the automaton numbered by $n$ and by $G_n$
the corresponding group of tree automorphisms. Sometimes, when the
context is clear, we use just the number to refer to the
corresponding automaton or group.

The following operations on automata change neither the group
generated by this automaton, nor, essentially, the action of the
group on the tree.
\begin{enumerate}
\item[($i$)] passing to inverses of all generators
\item[($ii$)] permuting the states of the automaton
\item[($iii$)] permuting the letters of the alphabet
\end{enumerate}

\begin{definition}
Two automata $\mathcal A$ and $\mathcal B$ that can be obtained from
one another using a composition of the operations ($i$)--($iii$),
are called \emph{symmetric}.
\end{definition}

\begin{definition}
If the minimization of an automaton $\mathcal A$ is symmetric to the
minimization of an automaton $\mathcal B$, we say that the automata
$\mathcal A$ and $\mathcal B$ are \emph{minimally symmetric} and
write $\mathcal A\sim\mathcal B$.
\end{definition}

Another equivalence relation we consider is the isomorphism of the
groups generated by the automata. The minimal symmetry relation is a
refinement of the isomorphism relation, since the same abstract
group may have different actions on the binary tree.

There are $194$ classes of automata, which are pairwise not
minimally symmetric, $10$ of which are minimally symmetric to
automata with fewer than $3$ states. These $10$ classes of automata
are subject of Theorem~\ref{thm_class2states}, which states that
they generate $6$ different groups.

At present, it is known that there are at most $124$ non-isomorphic
groups in the considered class.

\section{Selected groups from the class}\label{sec_tables}

In this section we provide information about selected groups in the
class of all groups generated by $(3,2)$ automata. The groups are
selected in such a way that the corresponding proofs in
Section~\ref{sec_proofs} show most of the main methods and ideas
that were used for the whole class.

The following notation is used:
\begin{itemize}
\item
Rels - this is a list of some relators in the group. All independent
relators up to length $20$ are included. On some situations
additional longer relators are included. For $G_{753}$ and $G_{858}$
there are no relators of length up to $20$ and the relators provided
in the table are not necessarily the shortest. In many cases, the
given relations are not sufficient (for example, some of the groups
are not finitely presented).
\item
SF - these numbers represent the size of the factors
$G/\Stg(n)$, for $n\geq0$.
\item Gr - these numbers represent the values of the growth
function $\gamma_G(n)$, for $n\geq0$, and generating system $a$,
$b$, $c$.
\end{itemize}

Finally, for each automaton in the list a histogram for the spectral
density of the operator $M_9$ acting on level $9$ of the tree is
shown.

In some cases, in order to show the main ways to prove the group
isomorphism, we provide several different automata generating the
same group. \vspace{0.5cm}

\noindent\begin{tabular}{@{}p{172pt}p{174pt}@{}}
\vspace{0.1cm}

\textbf{Automaton number $739$} \vspace{.2cm}

\begin{tabular}{@{}p{48pt}p{200pt}}

$a=\sigma(a,a)$

$b=(b,a)$

$c=(a,a)$& Group: $C_2\ltimes\bigl(C_2\wr\mathbb{Z}\bigr)$

Contracting: \emph{yes}

Self-replicating: \emph{no}\\
\end{tabular}

\begin{minipage}{230pt}
\setlength{\fboxrule}{0cm}
\begin{window}[5,r,{\fbox{\shortstack{\hspace{1.6cm}~\phantom{a}\\ \vspace{3.8cm}\\ \phantom{a}}}},{}]
Rels: $a^{2}$, $b^{2}$, $c^{2}$, $acac$, $acbabcabab$,
$acbabacbab$,\\
$abacbcacbcab$, $acbacbabcabc$, $acbcbabcabcbab$,\\
$acbcbabacbcbab$, $acbacbcacbcabc$, $abcbacbcacbcbcab$,\\
$acbcbacbabcbcabc$, $acbcbcbabcabcbcbab$,\\ $acbcbcbabacbcbcbab$,
$abcbacbcbabcbcabcb$,\\ $acbcbacbcacbcbcabc$
\\
SF: $2^0$,$2^{1}$,$2^{3}$,$2^{6}$,$2^{8}$,$2^{10}$,$2^{12}$,$2^{14}$,$2^{16}$\\
Gr: 1,4,9,17,30,47,68,93,122,155,192\\
\end{window}
\end{minipage}
& \hfill~

\hfill
\begin{picture}(1450,1090)(0,130)
\put(200,200){\circle{200}} 
\put(1200,200){\circle{200}}
\put(700,1070){\circle{200}}
\allinethickness{0.2mm} \put(45,280){$a$} \put(1280,280){$b$}
\put(820,1035){$c$}
\put(164,165){$\sigma$}  
\put(1164,152){$1$}       
\put(664,1022){$1$}       
\put(100,100){\arc{200}{0}{4.71}}     
\path(46,216)(100,200)(55,167)        
\spline(287,150)(700,0)(1113,150)     
\path(325,109)(287,150)(343,156)      
\put(1300,100){\arc{200}{4.71}{3.14}} 
\path(1345,167)(1300,200)(1354,216)     
\spline(650,983)(250,287)      
\path(297,318)(250,287)(253,343)      
\put(190,10){$_{0,1}$}  
\put(1155,10){$_0$}  
\put(680,77){$_1$}   
\put(455,585){$_{0,1}$}  
\end{picture}

\vspace{.3cm}

\hfill\epsfig{file=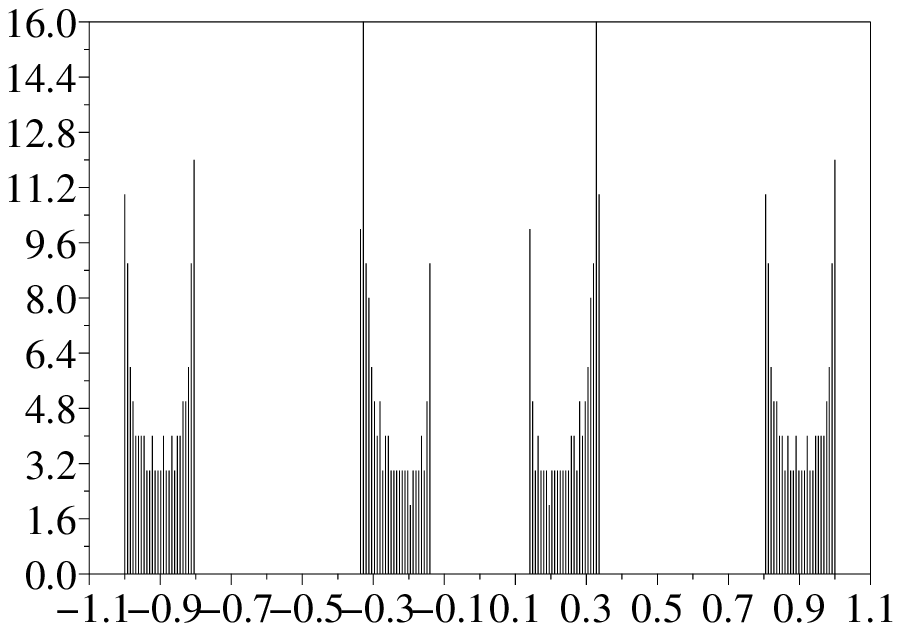,height=150pt}\\
\hline
\end{tabular}
\vspace{.3cm}

\noindent\begin{tabular}{@{}p{172pt}p{174pt}@{}} \textbf{Automaton
number $744$} \vspace{.2cm}

\begin{tabular}{@{}p{48pt}p{200pt}}

$a=\sigma(c,b)$

$b=(b,a)$

$c=(a,a)$& Group:

Contracting: \emph{no}

Self-replicating: \emph{yes}\\
\end{tabular}

\begin{minipage}{230pt}
\setlength{\fboxrule}{0cm}
\begin{window}[5,r,{\fbox{\shortstack{\hspace{1.6cm}~\phantom{a}\\ \vspace{3.8cm}\\ \phantom{a}}}},{}]
Rels:
$abcb^{-1}ac^{-1}a^{-2}bcb^{-1}ac^{-1}aca^{-1}bc^{-1}b^{-1}ca^{-1}bc^{-1}b^{-1}$,
$abcb^{-1}ac^{-1}a^{-2}bcb^{-1}ab^{-1}aca^{-1}bc^{-1}a^{-1}bc^{-1}b^{-1}$,
$abcb^{-1}ab^{-1}a^{-2}bcb^{-1}ac^{-1}aba^{-1}bc^{-1}b^{-1}ca^{-1}bc^{-1}b^{-1}$,
$abcb^{-1}ab^{-1}a^{-2}bcb^{-1}ab^{-1}aba^{-1}bc^{-1}a^{-1}bc^{-1}b^{-1}$
\\
SF: $2^0$,$2^{1}$,$2^{3}$,$2^{6}$,$2^{12}$,$2^{23}$,$2^{45}$,$2^{88}$,$2^{174}$\\
Gr: 1,7,37,187,937,4687\\
\end{window}
\end{minipage}
& \hfill~

\hfill
\begin{picture}(1450,1090)(0,130)
\put(200,200){\circle{200}} 
\put(1200,200){\circle{200}}
\put(700,1070){\circle{200}}
\allinethickness{0.2mm} \put(45,280){$a$} \put(1280,280){$b$}
\put(820,1035){$c$}
\put(164,165){$\sigma$}  
\put(1164,152){$1$}       
\put(664,1022){$1$}       
\put(300,200){\line(1,0){800}} 
\path(1050,225)(1100,200)(1050,175)   
\spline(200,300)(277,733)(613,1020)   
\path(559,1007)(613,1020)(591,969)    
\spline(287,150)(700,0)(1113,150)     
\path(325,109)(287,150)(343,156)      
\put(1300,100){\arc{200}{4.71}{3.14}} 
\path(1345,167)(1300,200)(1354,216)     
\spline(650,983)(250,287)      
\path(297,318)(250,287)(253,343)      
\put(230,700){$_0$} 
\put(680,240){$_1$} 
\put(1155,10){$_0$}  
\put(680,77){$_1$}   
\put(455,585){$_{0,1}$}  
\end{picture}

\vspace{.3cm}

\hfill\epsfig{file=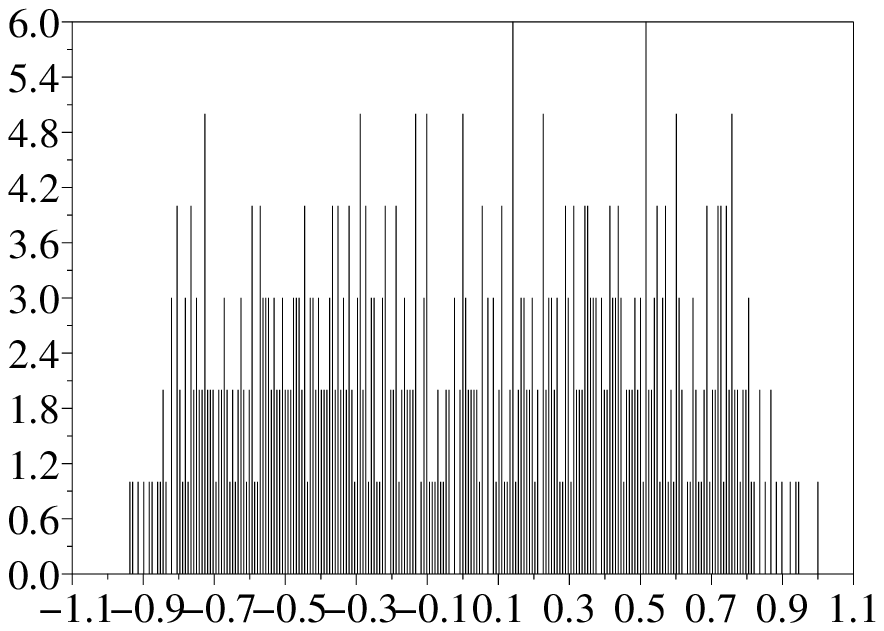,height=150pt}\\
\end{tabular}
\vspace{.3cm}

\noindent\begin{tabular}{@{}p{172pt}p{174pt}@{}} \textbf{Automaton
number $748$} \vspace{.2cm}

\begin{tabular}{@{}p{48pt}p{200pt}}

$a=\sigma(a,a)$

$b=(c,a)$

$c=(a,a)$& Group: $D_4\times C_2$

Contracting: \emph{yes}

Self-replicating: \emph{no}\\
\end{tabular}

\begin{minipage}{230pt}
\setlength{\fboxrule}{0cm}
\begin{window}[5,r,{\fbox{\shortstack{\hspace{1.6cm}~\phantom{a}\\ \vspace{3.8cm}\\ \phantom{a}}}},{}]
Rels: $a^{2}$, $b^{2}$, $c^{2}$, $acac$, $bcbc$, $abababab$
\\
SF: $2^0$,$2^{1}$,$2^{3}$,$2^{4}$,$2^{4}$,$2^{4}$,$2^{4}$,$2^{4}$,$2^{4}$\\
Gr: 1,4,8,12,15,16,16,16,16,16,16\\
\end{window}
\end{minipage}
& \hfill~

\hfill
\begin{picture}(1450,1090)(0,130)
\put(200,200){\circle{200}} 
\put(1200,200){\circle{200}}
\put(700,1070){\circle{200}}
\allinethickness{0.2mm} \put(45,280){$a$} \put(1280,280){$b$}
\put(820,1035){$c$}
\put(164,165){$\sigma$}  
\put(1164,152){$1$}       
\put(664,1022){$1$}       
\put(100,100){\arc{200}{0}{4.71}}     
\path(46,216)(100,200)(55,167)        
\spline(287,150)(700,0)(1113,150)     
\path(325,109)(287,150)(343,156)      
\spline(750,983)(1150,287)     
\path(753,927)(750,983)(797,952)      
\spline(650,983)(250,287)      
\path(297,318)(250,287)(253,343)      
\put(190,10){$_{0,1}$}  
\put(890,585){$_0$} 
\put(680,77){$_1$}   
\put(455,585){$_{0,1}$}  
\end{picture}

\vspace{.3cm}

\hfill\epsfig{file=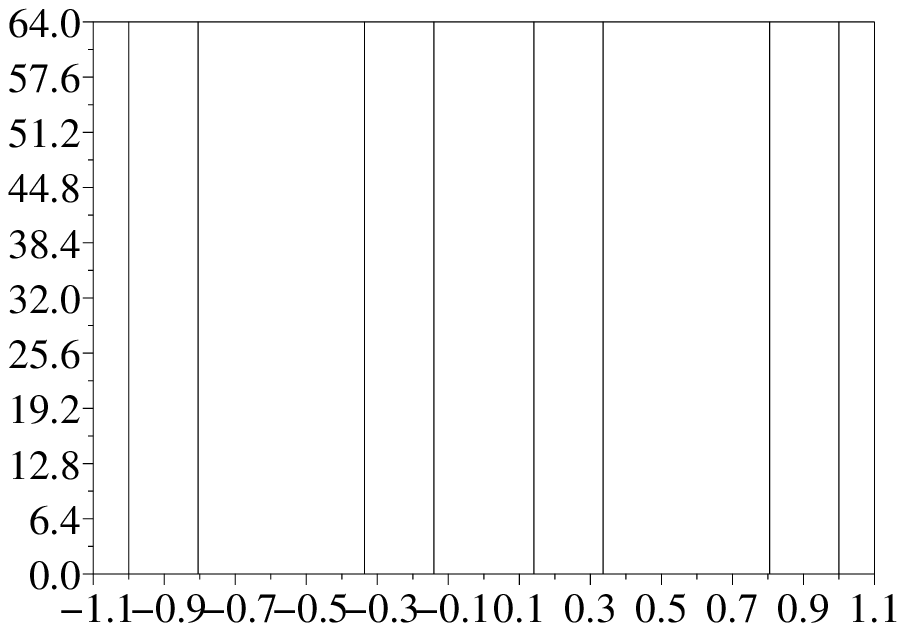,height=150pt}\\
\hline
\end{tabular}
\vspace{.3cm}

\noindent\begin{tabular}{@{}p{172pt}p{174pt}@{}} \textbf{Automaton
number $753$} \vspace{.2cm}

\begin{tabular}{@{}p{48pt}p{200pt}}

$a=\sigma(c,b)$

$b=(c,a)$

$c=(a,a)$& Group:

Contracting: \emph{no}

Self-replicating: \emph{yes}\\
\end{tabular}

\begin{minipage}{230pt}
\setlength{\fboxrule}{0cm}
\begin{window}[5,r,{\fbox{\shortstack{\hspace{1.6cm}~\phantom{a}\\ \vspace{3.8cm}\\ \phantom{a}}}},{}]
Rels:
$aba^{-1}b^{-1}ab^{-1}ca^{-1}ba^{-1}b^{-1}ab^{-1}cac^{-1}ba^{-1}bab^{-1}\cdot$\\
$a^{-1}c^{-1}ba^{-1}bab^{-1}$,
$aba^{-1}b^{-1}ab^{-1}ca^{-1}c^{-1}ba^{-1}c^{-1}b\cdot$\\
$ab^{-1}cac^{-1}ba^{-1}bab^{-1}a^{-1}c^{-1}ba^{-1}b^{-1}cab^{-1}c$,
$ac^{-1}ba^{-1}c^{-1}bab^{-1}ca^{-1}ba^{-1}b^{-1}ab^{-1}cac^{-1}ba^{-1}b^{-1}ca\cdot$\\
$b^{-1}ca^{-1}c^{-1}ba^{-1}bab^{-1}$
\\
SF: $2^0$,$2^{1}$,$2^{3}$,$2^{6}$,$2^{12}$,$2^{23}$,$2^{45}$,$2^{88}$,$2^{174}$\\
Gr: 1,7,37,187,937,4687\\
\end{window}
\end{minipage}
& \hfill~

\hfill
\begin{picture}(1450,1090)(0,130)
\put(200,200){\circle{200}} 
\put(1200,200){\circle{200}}
\put(700,1070){\circle{200}}
\allinethickness{0.2mm} \put(45,280){$a$} \put(1280,280){$b$}
\put(820,1035){$c$}
\put(164,165){$\sigma$}  
\put(1164,152){$1$}       
\put(664,1022){$1$}       
\put(300,200){\line(1,0){800}} 
\path(1050,225)(1100,200)(1050,175)   
\spline(200,300)(277,733)(613,1020)   
\path(559,1007)(613,1020)(591,969)    
\spline(287,150)(700,0)(1113,150)     
\path(325,109)(287,150)(343,156)      
\spline(750,983)(1150,287)     
\path(753,927)(750,983)(797,952)      
\spline(650,983)(250,287)      
\path(297,318)(250,287)(253,343)      
\put(230,700){$_0$} 
\put(680,240){$_1$} 
\put(890,585){$_0$} 
\put(680,77){$_1$}   
\put(455,585){$_{0,1}$}  
\end{picture}

\vspace{.3cm}

\hfill\epsfig{file=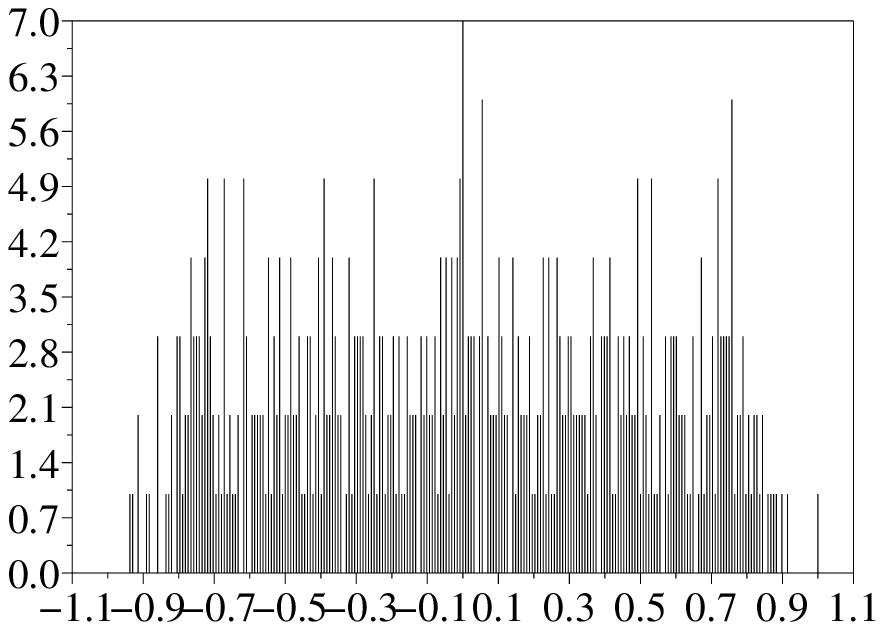,height=150pt}\\
\end{tabular}
\vspace{.3cm}

\noindent\begin{tabular}{@{}p{172pt}p{174pt}@{}} \textbf{Automaton
number $771$} \vspace{.2cm}

\begin{tabular}{@{}p{48pt}p{200pt}}

$a=\sigma(c,b)$

$b=(b,b)$

$c=(a,a)$& Group: $\mathbb Z^2$

Contracting: \emph{yes}

Self-replicating: \emph{yes}\\
\end{tabular}

\begin{minipage}{230pt}
\setlength{\fboxrule}{0cm}
\begin{window}[5,r,{\fbox{\shortstack{\hspace{1.6cm}~\phantom{a}\\ \vspace{3.8cm}\\ \phantom{a}}}},{}]
Rels: $b$, $a^{-1}c^{-1}ac$
\\
SF: $2^0$,$2^{1}$,$2^{2}$,$2^{3}$,$2^{4}$,$2^{5}$,$2^{6}$,$2^{7}$,$2^{8}$\\
Gr: 1,5,13,25,41,61,85,113,145,181,221\\
Limit space: $2$-dimensional torus $T_2$
\end{window}
\end{minipage}
& \hfill~

\hfill
\begin{picture}(1450,1090)(0,130)
\put(200,200){\circle{200}} 
\put(1200,200){\circle{200}}
\put(700,1070){\circle{200}}
\allinethickness{0.2mm} \put(45,280){$a$} \put(1280,280){$b$}
\put(820,1035){$c$}
\put(164,165){$\sigma$}  
\put(1164,152){$1$}       
\put(664,1022){$1$}       
\put(300,200){\line(1,0){800}} 
\path(1050,225)(1100,200)(1050,175)   
\spline(200,300)(277,733)(613,1020)   
\path(559,1007)(613,1020)(591,969)    
\put(1300,100){\arc{200}{4.71}{3.14}} 
\path(1345,167)(1300,200)(1354,216)     
\spline(650,983)(250,287)      
\path(297,318)(250,287)(253,343)      
\put(230,700){$_0$} 
\put(680,240){$_1$} 
\put(1080,10){$_{0,1}$}  
\put(455,585){$_{0,1}$}  
\end{picture}

\vspace{.3cm}

\hfill\epsfig{file=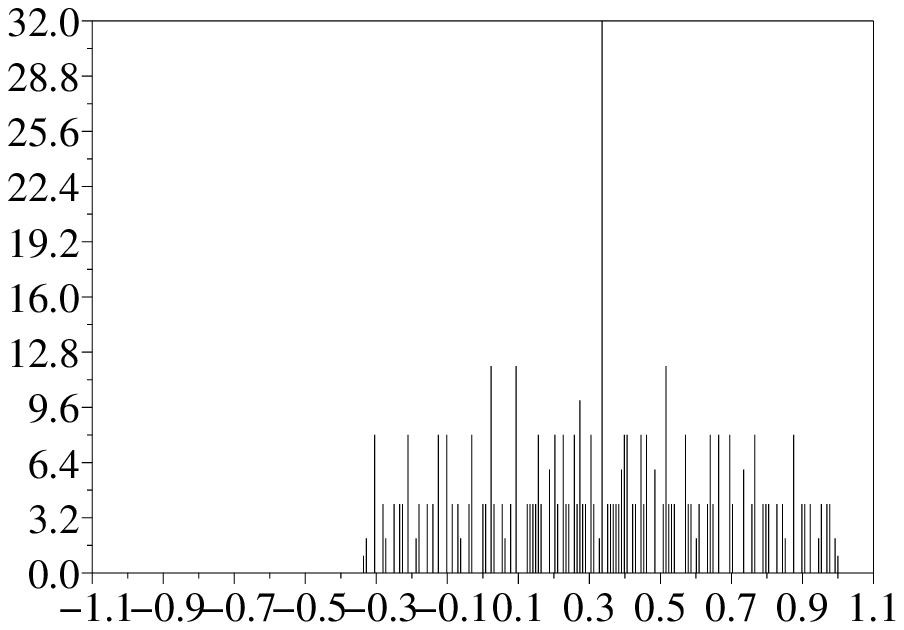,height=150pt}\\
\hline
\end{tabular}
\vspace{.3cm}

\noindent\begin{tabular}{@{}p{172pt}p{174pt}@{}} \textbf{Automata
number $775$ and $783$} \vspace{.2cm}

\begin{tabular}{@{}p{48pt}p{12pt}p{48pt}p{200pt}}

$a=\sigma(a,a)$

$b=(c,b)$

$c=(a,a)$&~

783:&

$a=\sigma(c,c)$

$b=(c,b)$

$c=(a,a)$&Group: $C_2\ltimes
IMG\left(\bigl(\frac{z-1}{z+1}\bigr)^2\right)$

Contracting: \emph{yes}

Self-replicating: \emph{yes}\\
\end{tabular}

\begin{minipage}{230pt}
\setlength{\fboxrule}{0cm}
\begin{window}[5,r,{\fbox{\shortstack{\hspace{1.6cm}~\phantom{a}\\ \vspace{3.8cm}\\ \phantom{a}}}},{}]
Rels: $a^{2}$, $b^{2}$, $c^{2}$, $acac$, $acbcbabcbcabcbabcb$,\\
$acbcbabcbacbcbabcb$, $abcbacbcbcacbcbcabcb$,\\
$acbcbacbcbabcbcabcbc$, $acbcbacbcbcacbcbcabcbc$
\\
SF: $2^0$,$2^{1}$,$2^{2}$,$2^{4}$,$2^{6}$,$2^{9}$,$2^{15}$,$2^{26}$,$2^{48}$\\
Gr: 1,4,9,17,30,51,85,140,229,367,579\\
Limit space:
\end{window}
\vspace{.1cm} \noindent\epsfig{file=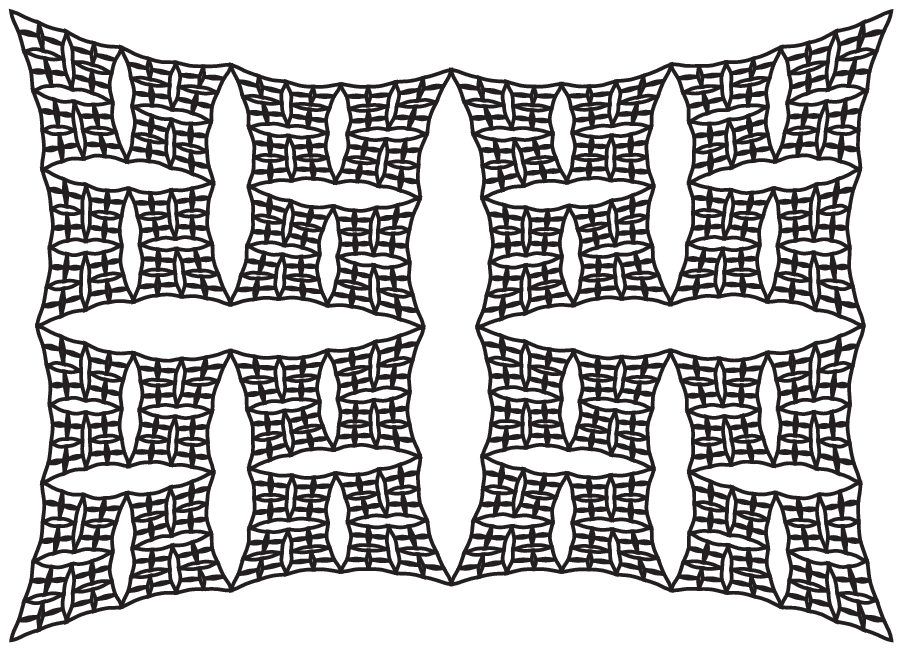,height=110pt}
\end{minipage} &
\hfill~

\hfill
\begin{picture}(1450,1090)(0,130)
\put(200,200){\circle{200}} 
\put(1200,200){\circle{200}}
\put(700,1070){\circle{200}}
\allinethickness{0.2mm} \put(45,280){$a$} \put(1280,280){$b$}
\put(820,1035){$c$}
\put(164,165){$\sigma$}  
\put(1164,152){$1$}       
\put(664,1022){$1$}       
\put(100,100){\arc{200}{0}{4.71}}     
\path(46,216)(100,200)(55,167)        
\put(1300,100){\arc{200}{4.71}{3.14}} 
\path(1345,167)(1300,200)(1354,216)     
\spline(750,983)(1150,287)     
\path(753,927)(750,983)(797,952)      
\spline(650,983)(250,287)      
\path(297,318)(250,287)(253,343)      
\put(190,10){$_{0,1}$}  
\put(890,585){$_0$} 
\put(1160,10){$_1$}  
\put(455,585){$_{0,1}$}  
\end{picture}

\vspace{.3cm}

\hfill\epsfig{file=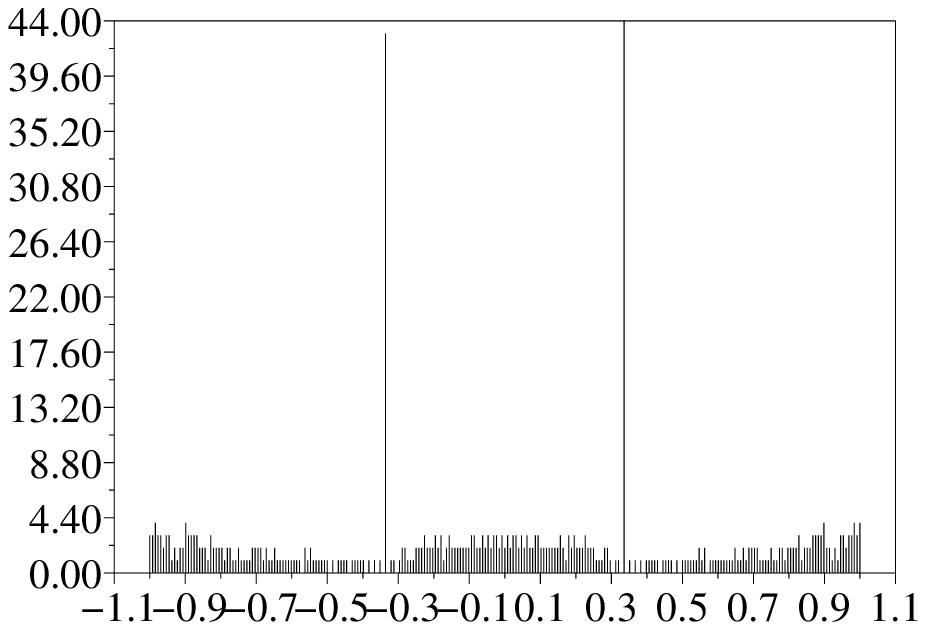,height=150pt}\\
\end{tabular}
\vspace{.3cm}

\noindent\begin{tabular}{@{}p{172pt}p{174pt}@{}} \textbf{Automaton
number $803$} \vspace{.2cm}

\begin{tabular}{@{}p{48pt}p{200pt}}

$a=\sigma(b,a)$

$b=(c,c)$

$c=(a,a)$& Group: $\mathbb Z^2$

Contracting: \emph{yes}

Self-replicating: \emph{yes}\\
\end{tabular}

\begin{minipage}{230pt}
\setlength{\fboxrule}{0cm}
\begin{window}[5,r,{\fbox{\shortstack{\hspace{1.6cm}~\phantom{a}\\ \vspace{3.8cm}\\ \phantom{a}}}},{}]
Rels: $cba^2$, $a^{-1}c^{-1}ac$
\\
SF: $2^0$,$2^{1}$,$2^{2}$,$2^{3}$,$2^{4}$,$2^{5}$,$2^{6}$,$2^{7}$,$2^{8}$\\
Gr: 1,7,21,43,73,111,157,211,273,343,421\\
Limit space: $2$-dimensional torus $T_2$
\end{window}
\end{minipage}
& \hfill~

\hfill
\begin{picture}(1450,1090)(0,130)
\put(200,200){\circle{200}} 
\put(1200,200){\circle{200}}
\put(700,1070){\circle{200}}
\allinethickness{0.2mm} \put(45,280){$a$} \put(1280,280){$b$}
\put(820,1035){$c$}
\put(164,165){$\sigma$}  
\put(1164,152){$1$}       
\put(664,1022){$1$}       
\put(100,100){\arc{200}{0}{4.71}}     
\path(46,216)(100,200)(55,167)        
\put(300,200){\line(1,0){800}} 
\path(1050,225)(1100,200)(1050,175)   
\spline(750,983)(1150,287)     
\path(753,927)(750,983)(797,952)      
\spline(650,983)(250,287)      
\path(297,318)(250,287)(253,343)      
\put(680,240){$_0$} 
\put(193,10){$_1$}  
\put(820,585){$_{0,1}$} 
\put(455,585){$_{0,1}$}  
\end{picture}

\vspace{.3cm}

\hfill\epsfig{file=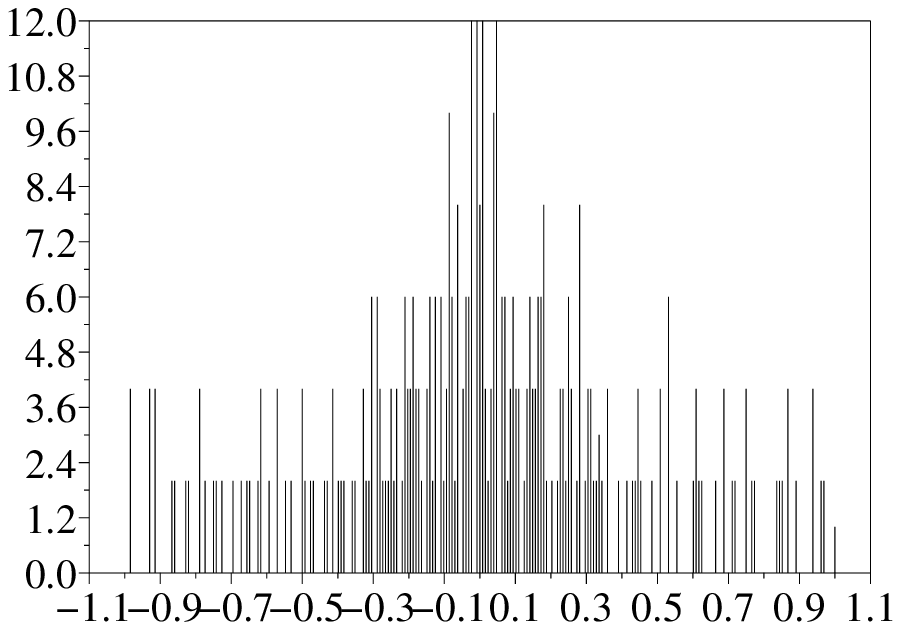,height=150pt}\\
\hline
\end{tabular}
\vspace{.3cm}

\noindent\begin{tabular}{@{}p{172pt}p{174pt}@{}} \textbf{Automaton
number $846$} \vspace{.2cm}

\begin{tabular}{@{}p{48pt}p{200pt}}

$a=\sigma(c,c)$

$b=(a,b)$

$c=(b,a)$& Group: $C_2\ast C_2\ast C_2$

Contracting: \emph{no}

Self-replicating: \emph{no}\\
\end{tabular}

\begin{minipage}{230pt}
\setlength{\fboxrule}{0cm}
\begin{window}[5,r,{\fbox{\shortstack{\hspace{1.6cm}~\phantom{a}\\ \vspace{3.8cm}\\ \phantom{a}}}},{}]
Rels: $a^{2}$, $b^{2}$, $c^{2}$
\\
SF: $2^0$,$2^{1}$,$2^{3}$,$2^{5}$,$2^{7}$,$2^{10}$,$2^{13}$,$2^{16}$,$2^{19}$\\
Gr: 1,4,10,22,46,94,190,382,766,1534\\
\end{window}
\end{minipage}
& \hfill~

\hfill
\begin{picture}(1450,1090)(0,130)
\put(200,200){\circle{200}} 
\put(1200,200){\circle{200}}
\put(700,1070){\circle{200}}
\allinethickness{0.2mm} \put(45,280){$a$} \put(1280,280){$b$}
\put(820,1035){$c$}
\put(164,165){$\sigma$}  
\put(1164,152){$1$}       
\put(664,1022){$1$}       
\spline(200,300)(277,733)(613,1020)   
\path(559,1007)(613,1020)(591,969)    
\spline(287,150)(700,0)(1113,150)     
\path(325,109)(287,150)(343,156)      
\put(1300,100){\arc{200}{4.71}{3.14}} 
\path(1345,167)(1300,200)(1354,216)     
\spline(650,983)(250,287)      
\path(297,318)(250,287)(253,343)      
\spline(1200,300)(1123,733)(787,1020) 
\path(1216,354)(1200,300)(1167,345)   
\put(150,700){$_{0,1}$} 
\put(680,77){$_0$}   
\put(1160,10){$_1$}  
\put(1115,700){$_0$}
\put(460,585){$_1$}  
\end{picture}

\vspace{.3cm}

\hfill\epsfig{file=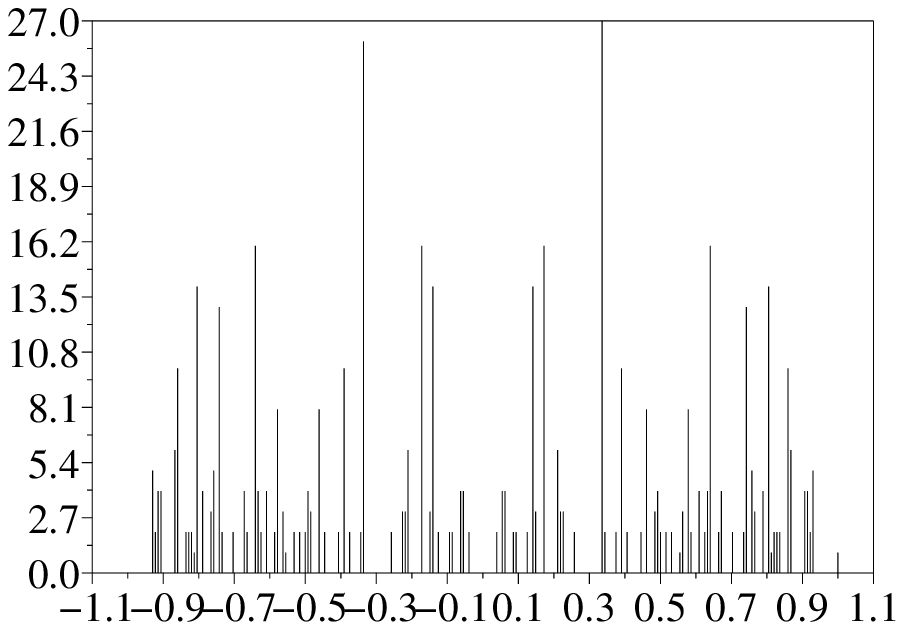,height=150pt}\\
\end{tabular}
\vspace{.3cm}

\noindent\begin{tabular}{@{}p{172pt}p{174pt}@{}} \textbf{Automaton
number $852$} \vspace{.2cm}

\begin{tabular}{@{}p{48pt}p{200pt}}

$a=\sigma(c,b)$

$b=(b,b)$

$c=(b,a)$& Group: \emph{Basilica Group} $\B=IMG(z^2-1)$

Contracting: \emph{yes}

Self-replicating: \emph{yes}\\
\end{tabular}

\begin{minipage}{230pt}
\setlength{\fboxrule}{0cm}
\begin{window}[5,r,{\fbox{\shortstack{\hspace{1.6cm}~\phantom{a}\\ \vspace{3.8cm}\\ \phantom{a}}}},{}]
Rels: $b$, $a^{-1}c^{-1}ac^{-1}a^{-1}cac$,
$a^{-1}c^{-2}ac^{-1}a^{-1}c^{2}ac$,\\
$a^{-1}c^{-1}ac^{-2}a^{-1}cac^{2}$,
$a^{-2}c^{-1}a^{-2}ca^{2}c^{-1}a^{2}c$,\\
$a^{-1}c^{-3}ac^{-1}a^{-1}c^{3}ac$,
$a^{-1}c^{-2}ac^{-2}a^{-1}c^{2}ac^{2}$,\\
$a^{-1}c^{-1}ac^{-3}a^{-1}cac^{3}$
\\
SF: $2^0$,$2^{1}$,$2^{3}$,$2^{6}$,$2^{12}$,$2^{23}$,$2^{45}$,$2^{88}$,$2^{174}$\\
Gr: 1,5,17,53,153,421,1125,2945,7545\\
Limit space:
\end{window}
\vspace{.5cm}

\noindent\epsfig{file=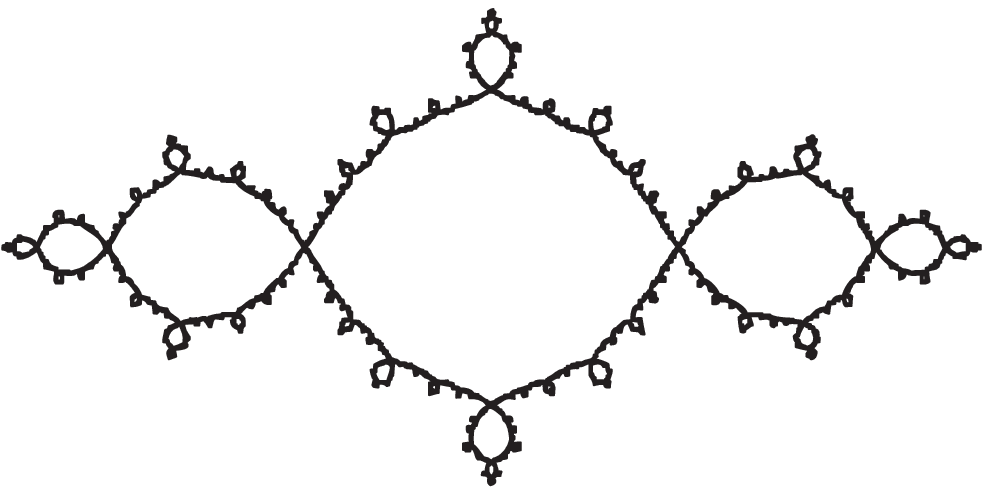,height=80pt}

\end{minipage}
& \hfill~

\hfill
\begin{picture}(1450,1090)(0,130)
\put(200,200){\circle{200}} 
\put(1200,200){\circle{200}}
\put(700,1070){\circle{200}}
\allinethickness{0.2mm} \put(45,280){$a$} \put(1280,280){$b$}
\put(820,1035){$c$}
\put(164,165){$\sigma$}  
\put(1164,152){$1$}       
\put(664,1022){$1$}       
\put(300,200){\line(1,0){800}} 
\path(1050,225)(1100,200)(1050,175)   
\spline(200,300)(277,733)(613,1020)   
\path(559,1007)(613,1020)(591,969)    
\put(1300,100){\arc{200}{4.71}{3.14}} 
\path(1345,167)(1300,200)(1354,216)     
\spline(650,983)(250,287)      
\path(297,318)(250,287)(253,343)      
\spline(1200,300)(1123,733)(787,1020) 
\path(1216,354)(1200,300)(1167,345)   
\put(230,700){$_0$} 
\put(680,240){$_1$} 
\put(1080,10){$_{0,1}$}  
\put(1115,700){$_0$}
\put(460,585){$_1$}  
\end{picture}

\vspace{.3cm}

\hfill\epsfig{file=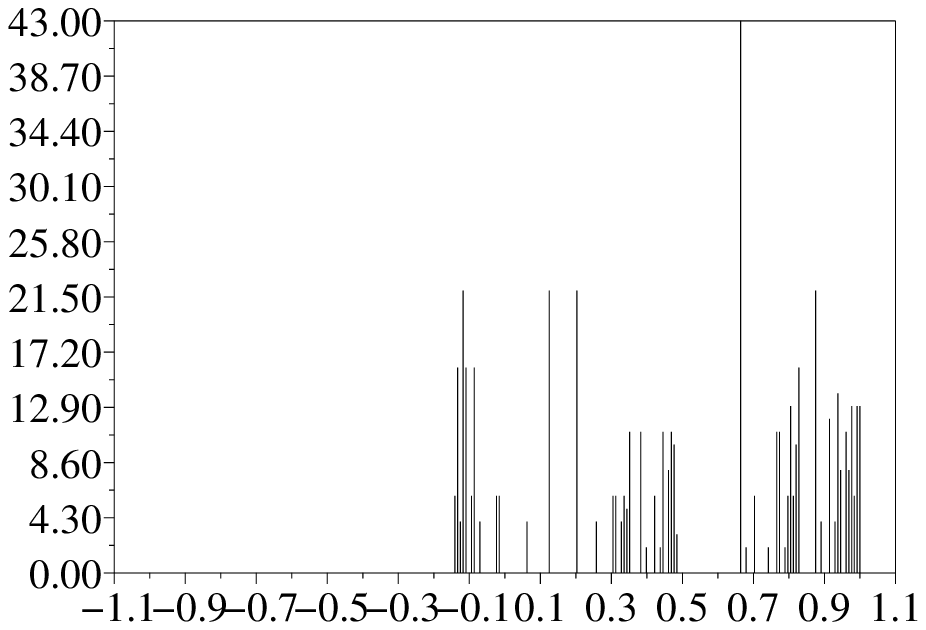,height=150pt}\\
\hline
\end{tabular}
\vspace{.3cm}

\noindent\begin{tabular}{@{}p{172pt}p{174pt}@{}} \textbf{Automaton
number $857$} \vspace{.2cm}

\begin{tabular}{@{}p{48pt}p{200pt}}

$a=\sigma(b,a)$

$b=(c,b)$

$c=(b,a)$& Group:

Contracting: \emph{no}

Self-replicating: \emph{yes}\\
\end{tabular}

\begin{minipage}{230pt}
\setlength{\fboxrule}{0cm}
\begin{window}[5,r,{\fbox{\shortstack{\hspace{1.6cm}~\phantom{a}\\ \vspace{3.8cm}\\ \phantom{a}}}},{}]
Rels: $a^{-1}ca^{-1}c$, $a^{-1}ba^{-1}ba^{-1}ba^{-1}b$,\\
$a^{-1}b^{-1}aca^{-1}b^{-1}ac$, $a^{-1}b^{-1}a^{2}c^{-1}b^{-1}ac$,
$a^{-2}bca^{-2}bc$, $b^{-1}cb^{-1}cb^{-1}cb^{-1}c$,
$a^{-1}ba^{-1}bc^{-1}ac^{-1}ba^{-1}b$,\\
$a^{-1}cac^{-1}b^{-1}aca^{-1}b^{-1}c$,
$a^{-1}bac^{-2}ac^{-1}bca^{-1}$
\\
SF: $2^0$,$2^{1}$,$2^{3}$,$2^{7}$,$2^{13}$,$2^{25}$,$2^{47}$,$2^{90}$,$2^{176}$\\
Gr: 1,7,35,165,758,3460\\
\end{window}
\end{minipage}
& \hfill~

\hfill
\begin{picture}(1450,1090)(0,130)
\put(200,200){\circle{200}} 
\put(1200,200){\circle{200}}
\put(700,1070){\circle{200}}
\allinethickness{0.2mm} \put(45,280){$a$} \put(1280,280){$b$}
\put(820,1035){$c$}
\put(164,165){$\sigma$}  
\put(1164,152){$1$}       
\put(664,1022){$1$}       
\put(100,100){\arc{200}{0}{4.71}}     
\path(46,216)(100,200)(55,167)        
\put(300,200){\line(1,0){800}} 
\path(1050,225)(1100,200)(1050,175)   
\put(1300,100){\arc{200}{4.71}{3.14}} 
\path(1345,167)(1300,200)(1354,216)     
\spline(750,983)(1150,287)     
\path(753,927)(750,983)(797,952)      
\spline(650,983)(250,287)      
\path(297,318)(250,287)(253,343)      
\spline(1200,300)(1123,733)(787,1020) 
\path(1216,354)(1200,300)(1167,345)   
\put(680,240){$_0$} 
\put(193,10){$_1$}  
\put(890,585){$_0$} 
\put(1160,10){$_1$}  
\put(1115,700){$_0$}
\put(460,585){$_1$}  
\end{picture}

\vspace{.3cm}

\hfill\epsfig{file=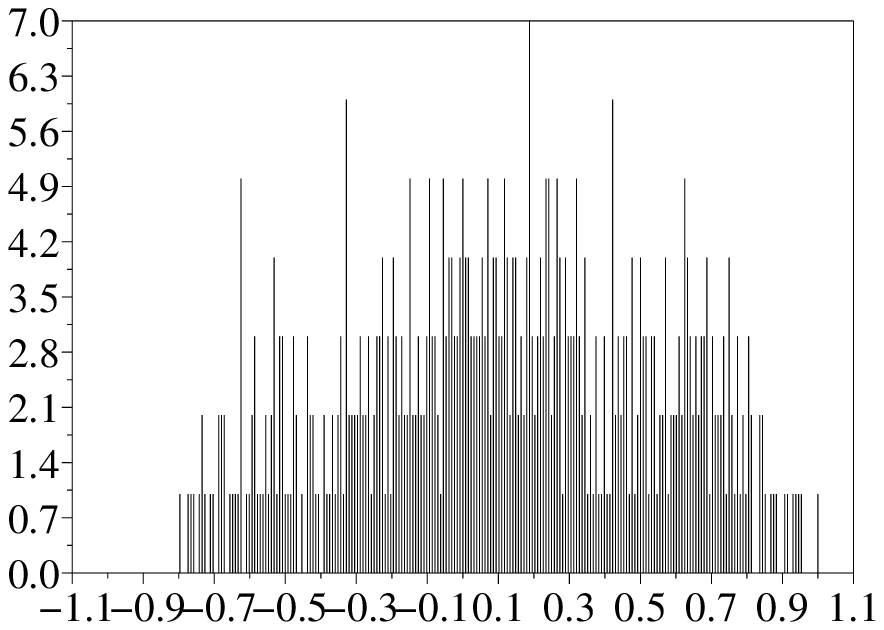,height=150pt}\\
\end{tabular}
\vspace{.3cm}

\noindent\begin{tabular}{@{}p{172pt}p{174pt}@{}} \textbf{Automaton
number $858$} \vspace{.2cm}

\begin{tabular}{@{}p{48pt}p{200pt}}

$a=\sigma(c,a)$

$b=(c,b)$

$c=(b,a)$& Group:

Contracting: \emph{no}

Self-replicating: \emph{yes}\\
\end{tabular}

\begin{minipage}{230pt}
\setlength{\fboxrule}{0cm}
\begin{window}[5,r,{\fbox{\shortstack{\hspace{1.6cm}~\phantom{a}\\ \vspace{3.8cm}\\ \phantom{a}}}},{}]
Rels:
$abca^{-1}c^{-1}ab^{-1}a^2c^{-1}b^{-1}a^{-1}bca^{-1}c^{-1}ab^{-1}a^2c^{-1}\cdot$\\
$b^{-1}abca^{-2}ba^{-1}cac^{-1}b^{-1}a^{-1}bca^{-2}ba^{-1}cac^{-1}b^{-1}$,
$abca^{-1}c^{-1}ab^{-1}a^2c^{-1}b^{-1}a^{-1}cba^{-1}b^{-1}ab^{-1}abca^{-2}b\cdot$\\
$a^{-1}cac^{-1}b^{-1}a^{-1}ba^{-1}bab^{-1}c^{-1}$
\\
SF: $2^0$,$2^{1}$,$2^{3}$,$2^{7}$,$2^{13}$,$2^{24}$,$2^{46}$,$2^{90}$,$2^{176}$\\
Gr: 1,7,37,187,937,4687\\
\end{window}
\end{minipage}
& \hfill~

\hfill
\begin{picture}(1450,1090)(0,130)
\put(200,200){\circle{200}} 
\put(1200,200){\circle{200}}
\put(700,1070){\circle{200}}
\allinethickness{0.2mm} \put(45,280){$a$} \put(1280,280){$b$}
\put(820,1035){$c$}
\put(164,165){$\sigma$}  
\put(1164,152){$1$}       
\put(664,1022){$1$}       
\put(100,100){\arc{200}{0}{4.71}}     
\path(46,216)(100,200)(55,167)        
\spline(200,300)(277,733)(613,1020)   
\path(559,1007)(613,1020)(591,969)    
\put(1300,100){\arc{200}{4.71}{3.14}} 
\path(1345,167)(1300,200)(1354,216)     
\spline(750,983)(1150,287)     
\path(753,927)(750,983)(797,952)      
\spline(650,983)(250,287)      
\path(297,318)(250,287)(253,343)      
\spline(1200,300)(1123,733)(787,1020) 
\path(1216,354)(1200,300)(1167,345)   
\put(230,700){$_0$} 
\put(193,10){$_1$}  
\put(890,585){$_0$} 
\put(1160,10){$_1$}  
\put(1115,700){$_0$}
\put(460,585){$_1$}  
\end{picture}

\vspace{.3cm}

\hfill\epsfig{file=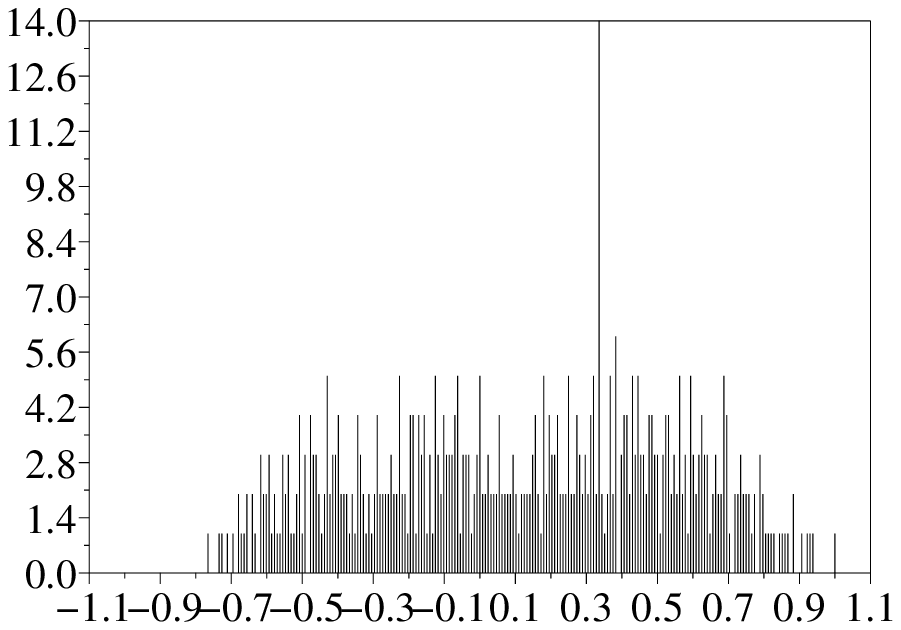,height=150pt}\\
\hline
\end{tabular}
\vspace{.3cm}

\noindent\begin{tabular}{@{}p{172pt}p{174pt}@{}} \textbf{Automaton
number $870$} \vspace{.2cm}

\begin{tabular}{@{}p{48pt}p{200pt}}

$a=\sigma(c,b)$

$b=(a,c)$

$c=(b,a)$& Group: \emph{Baumslag-Solitar group $BS(1,3)$}

Contracting: \emph{no}

Self-replicating: \emph{yes}\\
\end{tabular}

\begin{minipage}{230pt}
\setlength{\fboxrule}{0cm}
\begin{window}[5,r,{\fbox{\shortstack{\hspace{1.6cm}~\phantom{a}\\ \vspace{3.8cm}\\ \phantom{a}}}},{}]
Rels: $a^{-1}ca^{-1}b$, $(b^{-1}a)^b(b^{-1}a)^{-3}$
\\
SF: $2^0$,$2^{1}$,$2^{3}$,$2^{4}$,$2^{6}$,$2^{8}$,$2^{10}$,$2^{12}$,$2^{14}$\\
Gr: 1,7,33,127,433,1415\\
\end{window}
\end{minipage}
& \hfill~

\hfill
\begin{picture}(1450,1090)(0,130)
\put(200,200){\circle{200}} 
\put(1200,200){\circle{200}}
\put(700,1070){\circle{200}}
\allinethickness{0.2mm} \put(45,280){$a$} \put(1280,280){$b$}
\put(820,1035){$c$}
\put(164,165){$\sigma$}  
\put(1164,152){$1$}       
\put(664,1022){$1$}       
\put(300,200){\line(1,0){800}} 
\path(1050,225)(1100,200)(1050,175)   
\spline(200,300)(277,733)(613,1020)   
\path(559,1007)(613,1020)(591,969)    
\spline(287,150)(700,0)(1113,150)     
\path(325,109)(287,150)(343,156)      
\spline(750,983)(1150,287)     
\path(753,927)(750,983)(797,952)      
\spline(650,983)(250,287)      
\path(297,318)(250,287)(253,343)      
\spline(1200,300)(1123,733)(787,1020) 
\path(1216,354)(1200,300)(1167,345)   
\put(230,700){$_0$} 
\put(680,240){$_1$} 
\put(680,77){$_0$}   
\put(890,585){$_1$} 
\put(1115,700){$_0$}
\put(460,585){$_1$}  
\end{picture}

\vspace{.3cm}

\hfill\epsfig{file=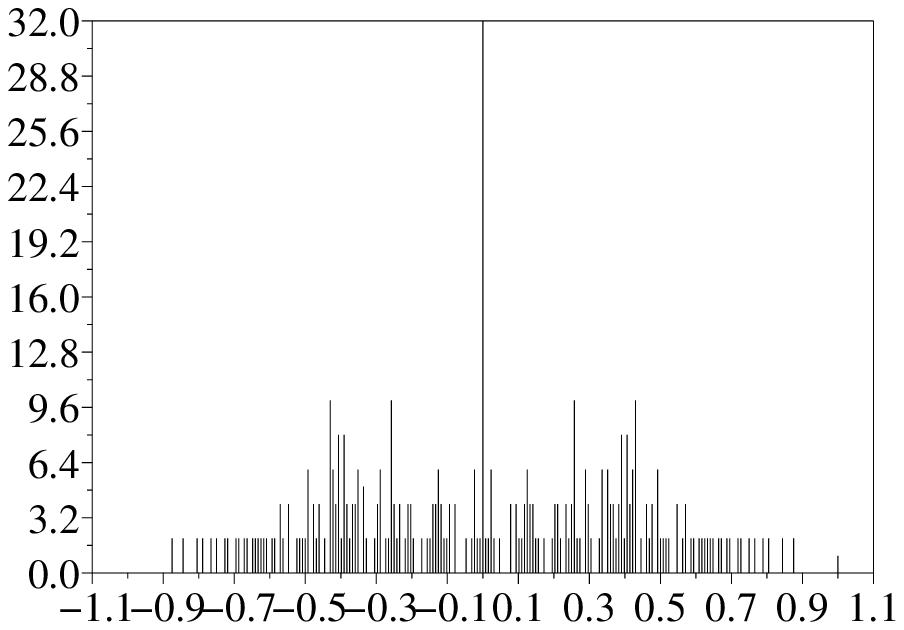,height=150pt}\\
\end{tabular}
\vspace{.3cm}

\noindent\begin{tabular}{@{}p{172pt}p{174pt}@{}} \textbf{Automaton
number $878$} \vspace{.2cm}

\begin{tabular}{@{}p{48pt}p{200pt}}

$a=\sigma(b,b)$

$b=(b,c)$

$c=(b,a)$& Group: $C_2\ltimes IMG(1-\frac1{z^2})$

Contracting: \emph{yes}

Self-replicating: \emph{yes}\\
\end{tabular}

\begin{minipage}{230pt}
\setlength{\fboxrule}{0cm}
\begin{window}[5,r,{\fbox{\shortstack{\hspace{1.6cm}~\phantom{a}\\ \vspace{3.8cm}\\ \phantom{a}}}},{}]
Rels: $a^{2}$, $b^{2}$, $c^{2}$, $abcabcacbacb$, $abcbcabcacbcbacb$
\\
SF: $2^0$,$2^{1}$,$2^{3}$,$2^{7}$,$2^{13}$,$2^{24}$,$2^{46}$,$2^{89}$,$2^{175}$\\
Gr: 1,4,10,22,46,94,184,352,664,1244,2296,4198,7612\\
Limit space:
\end{window}
\vspace{.5cm} \noindent\epsfig{file=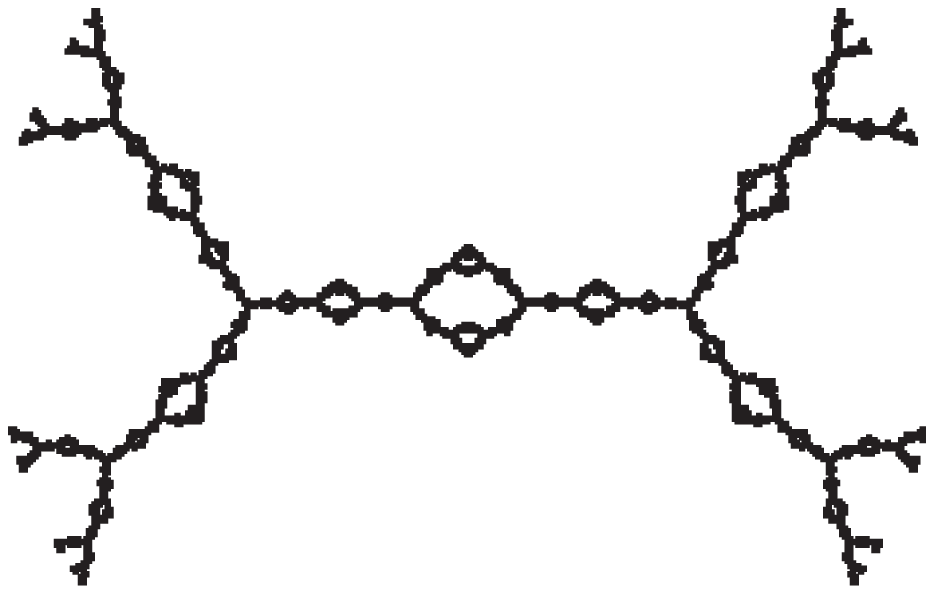,height=105pt}
\end{minipage} &
\hfill~

\hfill
\begin{picture}(1450,1090)(0,130)
\put(200,200){\circle{200}} 
\put(1200,200){\circle{200}}
\put(700,1070){\circle{200}}
\allinethickness{0.2mm} \put(45,280){$a$} \put(1280,280){$b$}
\put(820,1035){$c$}
\put(164,165){$\sigma$}  
\put(1164,152){$1$}       
\put(664,1022){$1$}       
\put(300,200){\line(1,0){800}} 
\path(1050,225)(1100,200)(1050,175)   
\put(1300,100){\arc{200}{4.71}{3.14}} 
\path(1345,167)(1300,200)(1354,216)     
\spline(750,983)(1150,287)     
\path(753,927)(750,983)(797,952)      
\spline(650,983)(250,287)      
\path(297,318)(250,287)(253,343)      
\spline(1200,300)(1123,733)(787,1020) 
\path(1216,354)(1200,300)(1167,345)   
\put(650,250){$_{0,1}$} 
\put(1155,10){$_0$}  
\put(890,585){$_1$} 
\put(1115,700){$_0$}
\put(460,585){$_1$}  
\end{picture}

\vspace{.3cm}

\hfill\epsfig{file=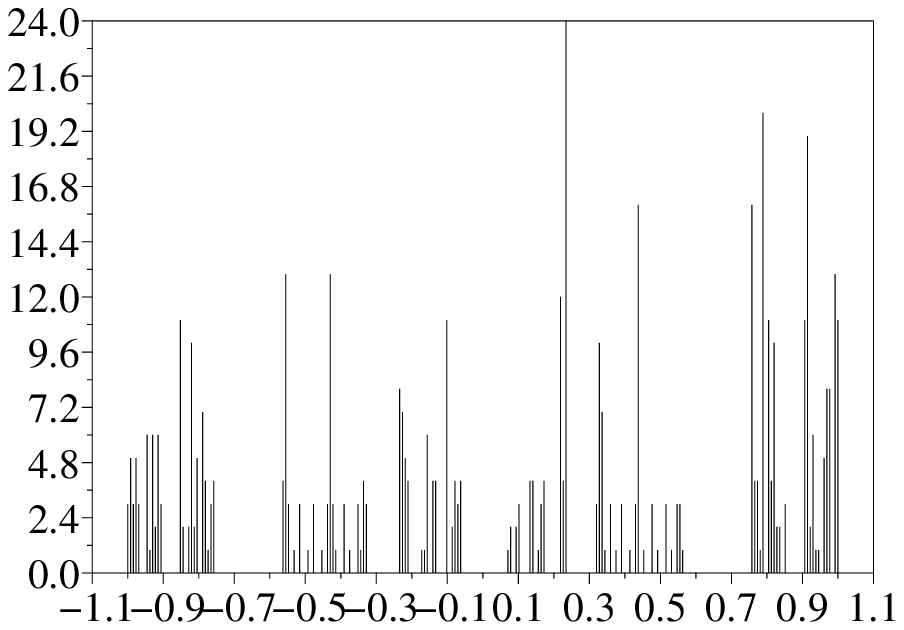,height=150pt}\\
\hline
\end{tabular}
\vspace{.3cm}

\noindent\begin{tabular}{@{}p{172pt}p{174pt}@{}} \textbf{Automaton
number $929$} \vspace{.2cm}

\begin{tabular}{@{}p{48pt}p{200pt}}

$a=\sigma(b,a)$

$b=(b,b)$

$c=(c,a)$& Group:

Contracting: \emph{no}

Self-replicating: \emph{yes}\\
\end{tabular}

\begin{minipage}{230pt}
\setlength{\fboxrule}{0cm}
\begin{window}[5,r,{\fbox{\shortstack{\hspace{1.6cm}~\phantom{a}\\ \vspace{3.8cm}\\ \phantom{a}}}},{}]
Rels: $b$, $a^{-3}cac^{-1}ac^{-1}ac$
\\
SF: $2^0$,$2^{1}$,$2^{3}$,$2^{6}$,$2^{12}$,$2^{23}$,$2^{45}$,$2^{88}$,$2^{174}$\\
Gr: 1,5,17,53,161,475,1387\\
\end{window}
\end{minipage}
& \hfill~

\hfill
\begin{picture}(1450,1090)(0,130)
\put(200,200){\circle{200}} 
\put(1200,200){\circle{200}}
\put(700,1070){\circle{200}}
\allinethickness{0.2mm} \put(45,280){$a$} \put(1280,280){$b$}
\put(820,1035){$c$}
\put(164,165){$\sigma$}  
\put(1164,152){$1$}       
\put(664,1022){$1$}       
\put(100,100){\arc{200}{0}{4.71}}     
\path(46,216)(100,200)(55,167)        
\put(300,200){\line(1,0){800}} 
\path(1050,225)(1100,200)(1050,175)   
\put(1300,100){\arc{200}{4.71}{3.14}} 
\path(1345,167)(1300,200)(1354,216)     
\spline(650,983)(250,287)      
\path(297,318)(250,287)(253,343)      
\put(700,1211){\arc{200}{2.36}{0.79}} 
\path(820,1168)(771,1141)(779,1196)   
\put(680,240){$_0$} 
\put(193,10){$_1$}  
\put(1080,10){$_{0,1}$}  
\put(545,1261){$_0$}  
\put(460,585){$_1$}  
\end{picture}

\vspace{.3cm}

\hfill\epsfig{file=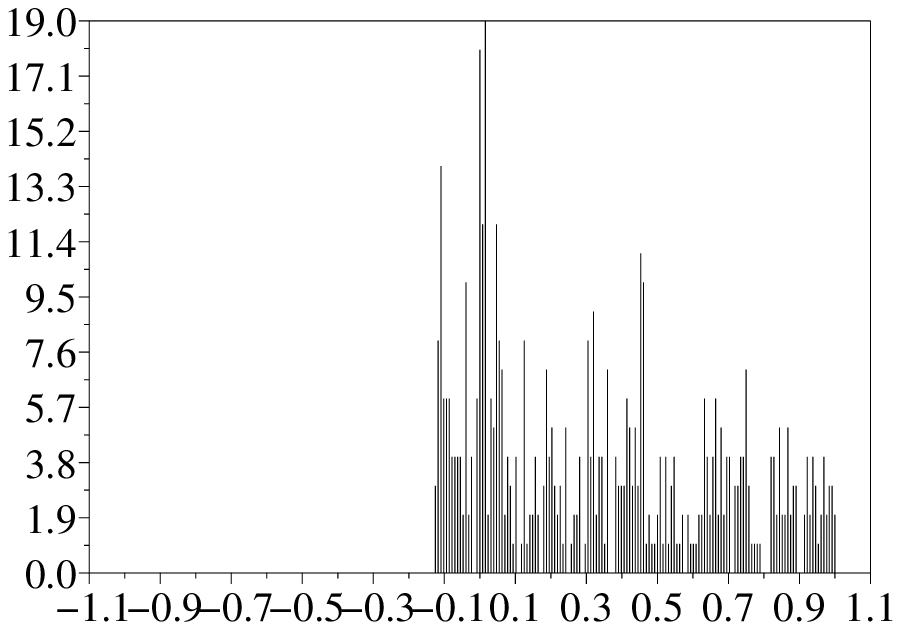,height=150pt}\\
\end{tabular}
\vspace{.3cm}

\noindent\begin{tabular}{@{}p{172pt}p{174pt}@{}} \textbf{Automaton
number $942$} \vspace{.2cm}

\begin{tabular}{@{}p{48pt}p{200pt}}

$a=\sigma(c,b)$

$b=(c,b)$

$c=(c,a)$& Group: Contains the lamplighter group

Contracting: \emph{no}

Self-replicating: \emph{yes}\\
\end{tabular}

\begin{minipage}{230pt}
\setlength{\fboxrule}{0cm}
\begin{window}[5,r,{\fbox{\shortstack{\hspace{1.6cm}~\phantom{a}\\ \vspace{3.8cm}\\ \phantom{a}}}},{}]
Rels: $a^{-1}ba^{-1}b$, $b^{-1}cb^{-1}c$, $b^{-1}ca^{-1}ba^{-1}c$,\\
$a^{-2}b^{2}a^{-1}b^{-1}ab$, $a^{-2}bab^{-2}ab$,
$b^{-2}c^{2}b^{-1}c^{-1}bc$,\\ $b^{-2}cbc^{-2}bc$,
$a^{-1}ca^{-1}ca^{-1}ca^{-1}c$, $b^{-1}ca^{-1}cb^{-1}ca^{-1}c$,\\
$a^{-1}bab^{-1}a^{-1}b^{2}a^{-1}$,
$b^{-1}cbc^{-1}b^{-1}c^{2}b^{-1}$,\\
$a^{-2}bac^{-1}bc^{-1}b^{-1}ab$,
$a^{-1}bc^{-1}ac^{-1}ac^{-1}ac^{-1}b$,\\
$b^{-1}ca^{-2}b^{2}a^{-1}b^{-1}ac$, $b^{-1}ca^{-2}bab^{-2}ac$,\\
$b^{-2}c^{2}a^{-1}ba^{-1}c^{-1}bc$, $b^{-1}cab^{-2}aba^{-2}c$,\\
$b^{-1}cab^{-1}a^{-1}b^{2}a^{-2}c$,\\
$a^{-1}bab^{-1}c^{-1}bc^{-1}aba^{-1}$,\\
$b^{-1}cbc^{-1}a^{-1}ba^{-1}c^{2}b^{-1}$
\\
SF: $2^0$,$2^{1}$,$2^{3}$,$2^{7}$,$2^{13}$,$2^{25}$,$2^{47}$,$2^{90}$,$2^{176}$\\
Gr: 1,7,33,143,597,2465\\
\end{window}
\end{minipage}
& \hfill~

\hfill
\begin{picture}(1450,1090)(0,130)
\put(200,200){\circle{200}} 
\put(1200,200){\circle{200}}
\put(700,1070){\circle{200}}
\allinethickness{0.2mm} \put(45,280){$a$} \put(1280,280){$b$}
\put(820,1035){$c$}
\put(164,165){$\sigma$}  
\put(1164,152){$1$}       
\put(664,1022){$1$}       
\put(300,200){\line(1,0){800}} 
\path(1050,225)(1100,200)(1050,175)   
\spline(200,300)(277,733)(613,1020)   
\path(559,1007)(613,1020)(591,969)    
\put(1300,100){\arc{200}{4.71}{3.14}} 
\path(1345,167)(1300,200)(1354,216)     
\spline(750,983)(1150,287)     
\path(753,927)(750,983)(797,952)      
\spline(650,983)(250,287)      
\path(297,318)(250,287)(253,343)      
\put(700,1211){\arc{200}{2.36}{0.79}} 
\path(820,1168)(771,1141)(779,1196)   
\put(230,700){$_0$} 
\put(680,240){$_1$} 
\put(890,585){$_0$} 
\put(1160,10){$_1$}  
\put(545,1261){$_0$}  
\put(460,585){$_1$}  
\end{picture}

\vspace{.3cm}

\hfill\epsfig{file=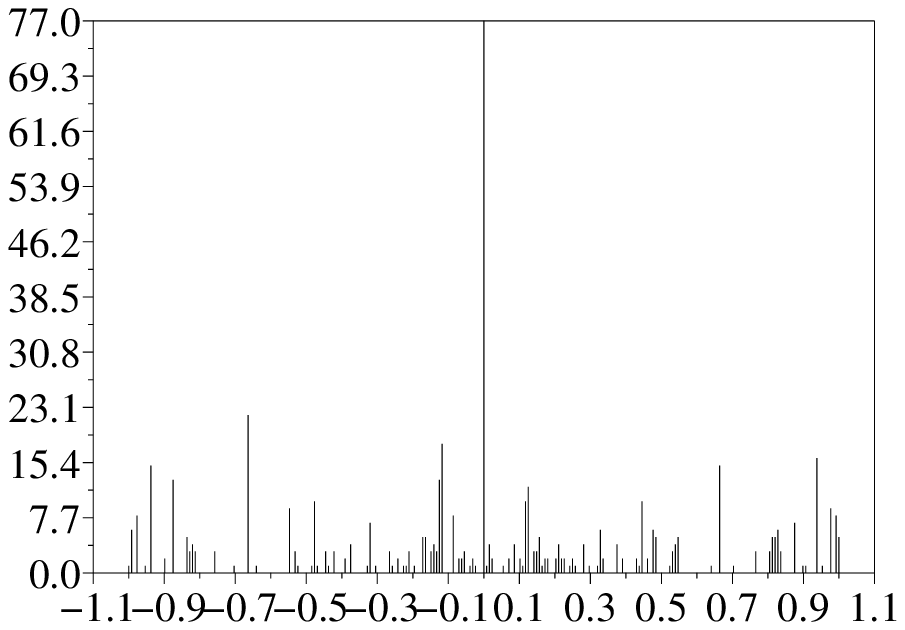,height=150pt}\\
\hline
\end{tabular}
\vspace{.3cm}

\noindent\begin{tabular}{@{}p{172pt}p{174pt}@{}} \textbf{Automaton
number $968$} \vspace{.2cm}

\begin{tabular}{@{}p{48pt}p{200pt}}

$a=\sigma(b,b)$

$b=(c,c)$

$c=(c,a)$& Group: contains $\Z^5$ as a subgroup of index $16$

Contracting: \emph{yes}

Self-replicating: \emph{no}\\
\end{tabular}

\begin{minipage}{230pt}
\setlength{\fboxrule}{0cm}
\begin{window}[5,r,{\fbox{\shortstack{\hspace{1.6cm}~\phantom{a}\\ \vspace{3.8cm}\\ \phantom{a}}}},{}]
Rels: $a^{2}$, $b^{2}$, $c^{2}$, $abcabcacbacb$, $acbcbabacbcbab$,\\
$acacbcbacacbcb$, $abcbcabcacbcbacb$, $acabacbabacabacbab$,\\
$acbabacacacbacacab$, $acacacacbacacacacb$,\\ $acbcbabcbacbcbabcb$,
$acbcacbcbacbcacbcb$
\\
SF: $2^0$,$2^{1}$,$2^{3}$,$2^{6}$,$2^{9}$,$2^{13}$,$2^{17}$,$2^{21}$,$2^{25}$\\
Gr: 1,4,10,22,46,94,184,338,600,1022\\
\end{window}
\end{minipage}
& \hfill~

\hfill
\begin{picture}(1450,1090)(0,130)
\put(200,200){\circle{200}} 
\put(1200,200){\circle{200}}
\put(700,1070){\circle{200}}
\allinethickness{0.2mm} \put(45,280){$a$} \put(1280,280){$b$}
\put(820,1035){$c$}
\put(164,165){$\sigma$}  
\put(1164,152){$1$}       
\put(664,1022){$1$}       
\put(300,200){\line(1,0){800}} 
\path(1050,225)(1100,200)(1050,175)   
\spline(750,983)(1150,287)     
\path(753,927)(750,983)(797,952)      
\spline(650,983)(250,287)      
\path(297,318)(250,287)(253,343)      
\put(700,1211){\arc{200}{2.36}{0.79}} 
\path(820,1168)(771,1141)(779,1196)   
\put(650,250){$_{0,1}$} 
\put(820,585){$_{0,1}$} 
\put(545,1261){$_0$}  
\put(460,585){$_1$}  
\end{picture}

\vspace{.3cm}

\hfill\epsfig{file=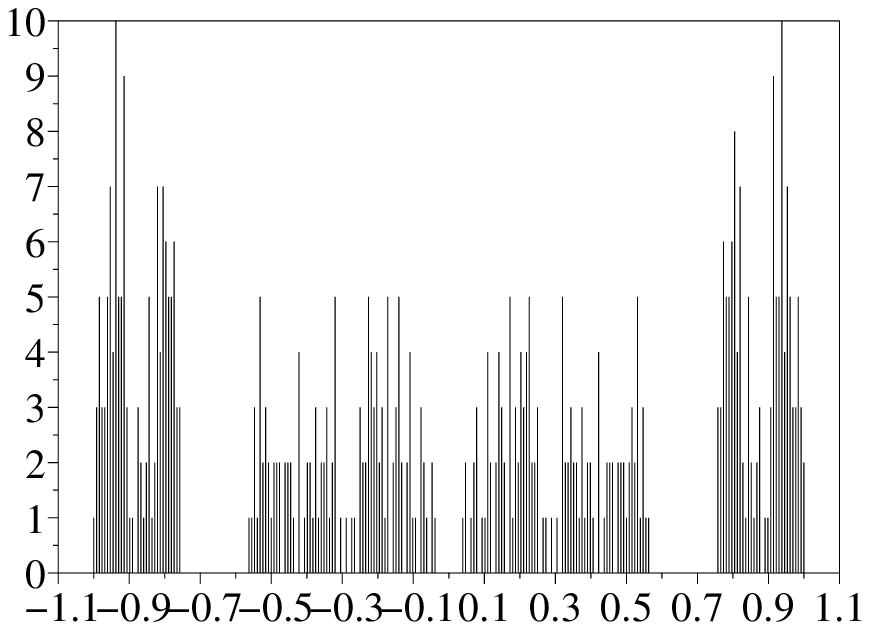,height=150pt}\\
\end{tabular}
\vspace{.3cm}

\noindent\begin{tabular}{@{}p{172pt}p{174pt}@{}} \textbf{Automaton
number $2205$} \vspace{.2cm}

\begin{tabular}{@{}p{48pt}p{200pt}}

$a=\sigma(c,c)$

$b=\sigma(b,a)$

$c=(a,a)$& Group: $C_2\ltimes
IMG\left(\bigl(\frac{z-1}{z+1}\bigr)^2\right)$

Contracting: \emph{yes}

Self-replicating: \emph{yes}\\
\end{tabular}

\begin{minipage}{230pt}
\setlength{\fboxrule}{0cm}
\begin{window}[5,r,{\fbox{\shortstack{\hspace{1.6cm}~\phantom{a}\\ \vspace{3.8cm}\\ \phantom{a}}}},{}]
Rels: $a^{2}$, $c^{2}$, $acac$, $acbcab$, $acbacb$,\\ $abcabc$,
$acb^{2}cab^{2}$,  $acb^{2}acb^{2}$, $ab^{2}cab^{2}c$,\\
$acb^{3}cab^{3}$, $acb^{3}acb^{3}$, $ab^{3}cab^{3}c$,
  $acb^{4}cab^{4}$,\\ $acb^{4}acb^{4}$, $ab^{4}cab^{4}c$\\
SF: $2^0$,$2^{1}$,$2^{2}$,$2^{4}$,$2^{6}$,$2^{9}$,$2^{15}$,$2^{26}$,$2^{48}$,$2^{91}$\\
Gr: 1,5,16,40,88,184,376,746,1458\\
Limit space:
\end{window}
\vspace{.1cm}
\noindent\epsfig{file=schmetterling1.eps,height=110pt}\end{minipage}
& \hfill~

\hfill
\begin{picture}(1450,1090)(0,130)
\put(200,200){\circle{200}} 
\put(1200,200){\circle{200}}
\put(700,1070){\circle{200}}
\allinethickness{0.2mm} \put(45,280){$a$} \put(1280,280){$b$}
\put(820,1035){$c$}
\put(164,165){$\sigma$}  
\put(1164,165){$\sigma$}  
\put(664,1022){$1$}       
\spline(200,300)(277,733)(613,1020)   
\path(559,1007)(613,1020)(591,969)    
\spline(287,150)(700,0)(1113,150)     
\path(325,109)(287,150)(343,156)      
\put(1300,100){\arc{200}{4.71}{3.14}} 
\path(1345,167)(1300,200)(1354,216)     
\spline(650,983)(250,287)      
\path(297,318)(250,287)(253,343)      
\put(150,700){$_{0,1}$} 
\put(1155,10){$_0$}  
\put(680,77){$_1$}   
\put(455,585){$_{0,1}$}  
\end{picture}

\vspace{.3cm}

\hfill\epsfig{file=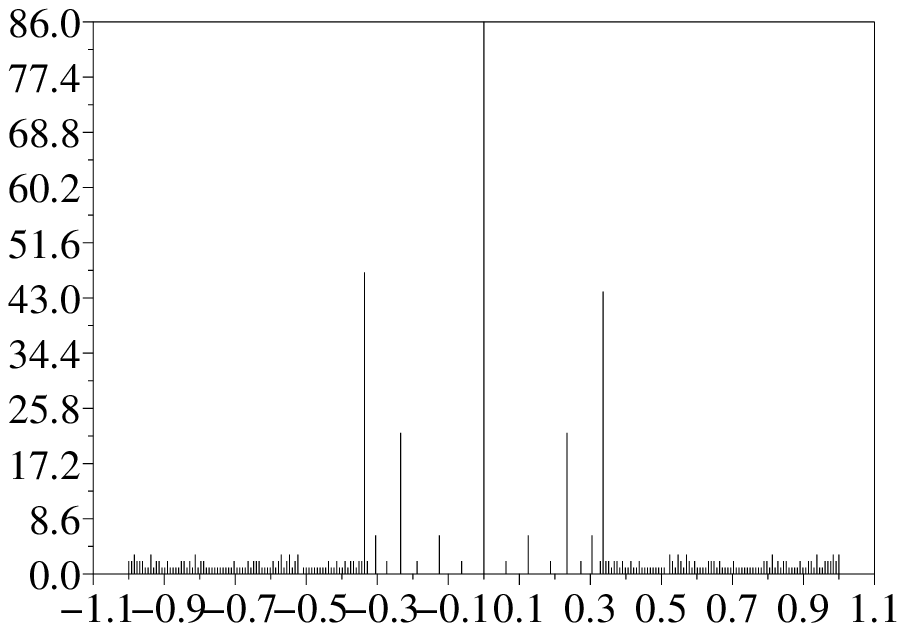,height=150pt}\\
\hline
\end{tabular}
\vspace{.3cm}

\noindent\begin{tabular}{@{}p{172pt}p{174pt}@{}} \textbf{Automaton
number $2212$} \vspace{.2cm}

\begin{tabular}{@{}p{48pt}p{200pt}}

$a=\sigma(a,c)$

$b=\sigma(c,a)$

$c=(a,a)$& Group: \emph{Klein bottle group}, virtually $\Z^2$

Contracting: \emph{yes}

Self-replicating: \emph{no}\\
\end{tabular}

\begin{minipage}{230pt}
\setlength{\fboxrule}{0cm}
\begin{window}[5,r,{\fbox{\shortstack{\hspace{1.6cm}~\phantom{a}\\ \vspace{3.8cm}\\ \phantom{a}}}},{}]
Rels: $ca^{2}$, $cb^{2}$
\\
SF: $2^0$,$2^{1}$,$2^{2}$,$2^{4}$,$2^{6}$,$2^{8}$,$2^{10}$,$2^{12}$,$2^{14}$\\
Gr: 1,7,19,37,61,91,127,169,217,271,331\\
\end{window}
\end{minipage}
& \hfill~

\hfill
\begin{picture}(1450,1090)(0,130)
\put(200,200){\circle{200}} 
\put(1200,200){\circle{200}}
\put(700,1070){\circle{200}}
\allinethickness{0.2mm} \put(45,280){$a$} \put(1280,280){$b$}
\put(820,1035){$c$}
\put(164,165){$\sigma$}  
\put(1164,165){$\sigma$}  
\put(664,1022){$1$}       
\put(100,100){\arc{200}{0}{4.71}}     
\path(46,216)(100,200)(55,167)        
\spline(200,300)(277,733)(613,1020)   
\path(559,1007)(613,1020)(591,969)    
\spline(287,150)(700,0)(1113,150)     
\path(325,109)(287,150)(343,156)      
\spline(750,983)(1150,287)     
\path(753,927)(750,983)(797,952)      
\spline(650,983)(250,287)      
\path(297,318)(250,287)(253,343)      
\put(190,10){$_0$}  
\put(230,700){$_1$} 
\put(890,585){$_0$} 
\put(680,77){$_1$}   
\put(455,585){$_{0,1}$}  
\end{picture}

\vspace{.3cm}

\hfill\epsfig{file=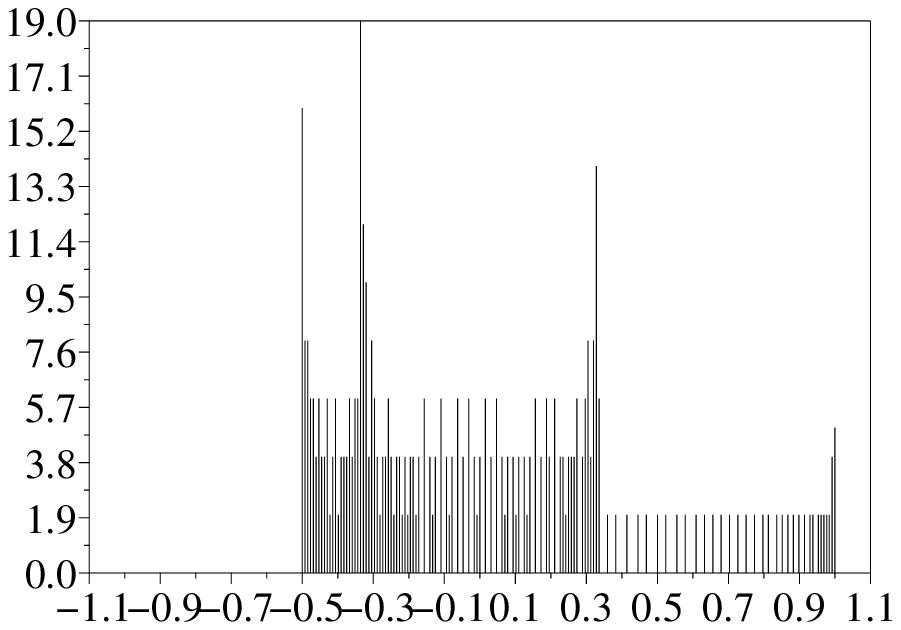,height=150pt}\\
\end{tabular}
\vspace{.3cm}

\noindent\begin{tabular}{@{}p{172pt}p{174pt}@{}} \textbf{Automaton
number $2240$} \vspace{.2cm}

\begin{tabular}{@{}p{48pt}p{200pt}}

$a=\sigma(b,c)$

$b=\sigma(c,b)$

$c=(a,a)$& Group: $F_3$~--- free of rank $3$ (Aleshin-Vorobets)

Contracting: \emph{no}

Self-replicating: \emph{no}\\
\end{tabular}

\begin{minipage}{230pt}
\setlength{\fboxrule}{0cm}
\begin{window}[5,r,{\fbox{\shortstack{\hspace{1.6cm}~\phantom{a}\\ \vspace{3.8cm}\\ \phantom{a}}}},{}]
Rels:\\
SF: $2^0$,$2^{1}$,$2^{2}$,$2^{4}$,$2^{7}$,$2^{10}$,$2^{14}$,$2^{21}$,$2^{34}$\\
Gr: 1,7,37,187,937,4687\\
\end{window}
\end{minipage}
& \hfill~

\hfill
\begin{picture}(1450,1090)(0,130)
\put(200,200){\circle{200}} 
\put(1200,200){\circle{200}}
\put(700,1070){\circle{200}}
\allinethickness{0.2mm} \put(45,280){$a$} \put(1280,280){$b$}
\put(820,1035){$c$}
\put(164,165){$\sigma$}  
\put(1164,165){$\sigma$}  
\put(664,1022){$1$}       
\put(300,200){\line(1,0){800}} 
\path(1050,225)(1100,200)(1050,175)   
\spline(200,300)(277,733)(613,1020)   
\path(559,1007)(613,1020)(591,969)    
\put(1300,100){\arc{200}{4.71}{3.14}} 
\path(1345,167)(1300,200)(1354,216)     
\spline(750,983)(1150,287)     
\path(753,927)(750,983)(797,952)      
\spline(650,983)(250,287)      
\path(297,318)(250,287)(253,343)      
\put(680,240){$_0$} 
\put(230,700){$_1$} 
\put(890,585){$_0$} 
\put(1160,10){$_1$}  
\put(455,585){$_{0,1}$}  
\end{picture}

\vspace{.3cm}

\hfill\epsfig{file=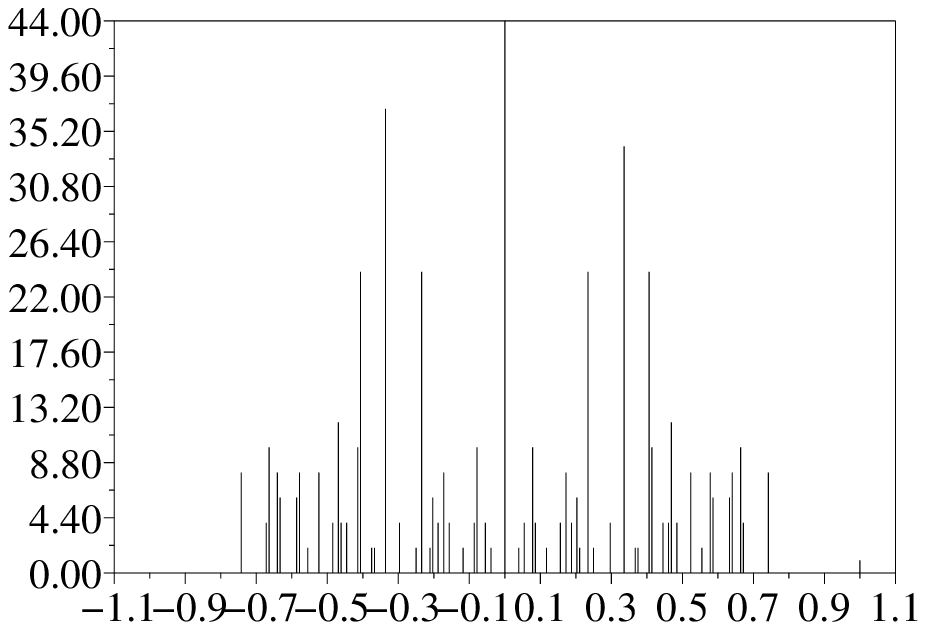,height=150pt}\\
\hline
\end{tabular}
\vspace{.3cm}

\noindent\begin{tabular}{@{}p{172pt}p{174pt}@{}} \textbf{Automaton
number $2277$} \vspace{.2cm}

\begin{tabular}{@{}p{48pt}p{200pt}}

$a=\sigma(c,c)$

$b=\sigma(a,a)$

$c=(b,a)$& Group: $C_2\ltimes(\mathbb Z\times\mathbb Z)$

Contracting: \emph{yes}

Self-replicating: \emph{yes}\\
\end{tabular}

\begin{minipage}{230pt}
\setlength{\fboxrule}{0cm}
\begin{window}[5,r,{\fbox{\shortstack{\hspace{1.6cm}~\phantom{a}\\ \vspace{3.8cm}\\ \phantom{a}}}},{}]
Rels: $a^{2}$, $b^{2}$, $c^{2}$, $acbacb$, $acbcbacbcb$,\\
$bacacbacac$, $acbcbcbacbcbcb$, $bacacacbacacac$,\\
$acbcbcbcbacbcbcbcb$, $bacacacacbacacacac$,\\
$acbcbcbcbcbacbcbcbcbcb$, $bacacacacacbacacacacac$,\\
$acbcbcbcbcbcbacbcbcbcbcbcb$
\\
SF: $2^0$,$2^{1}$,$2^{2}$,$2^{4}$,$2^{5}$,$2^{6}$,$2^{7}$,$2^{8}$,$2^{9}$\\
Gr: 1,4,10,19,31,46,64,85,109,136,166\\
Limit space: $2$-dimensional sphere $S_2$
\end{window}
\end{minipage}
& \hfill~

\hfill
\begin{picture}(1450,1090)(0,130)
\put(200,200){\circle{200}} 
\put(1200,200){\circle{200}}
\put(700,1070){\circle{200}}
\allinethickness{0.2mm} \put(45,280){$a$} \put(1280,280){$b$}
\put(820,1035){$c$}
\put(164,165){$\sigma$}  
\put(1164,165){$\sigma$}  
\put(664,1022){$1$}       
\spline(200,300)(277,733)(613,1020)   
\path(559,1007)(613,1020)(591,969)    
\spline(287,150)(700,0)(1113,150)     
\path(325,109)(287,150)(343,156)      
\spline(650,983)(250,287)      
\path(297,318)(250,287)(253,343)      
\spline(1200,300)(1123,733)(787,1020) 
\path(1216,354)(1200,300)(1167,345)   
\put(150,700){$_{0,1}$} 
\put(650,87){$_{0,1}$}   
\put(1115,700){$_0$}
\put(460,585){$_1$}  
\end{picture}

\vspace{.3cm}

\hfill\epsfig{file=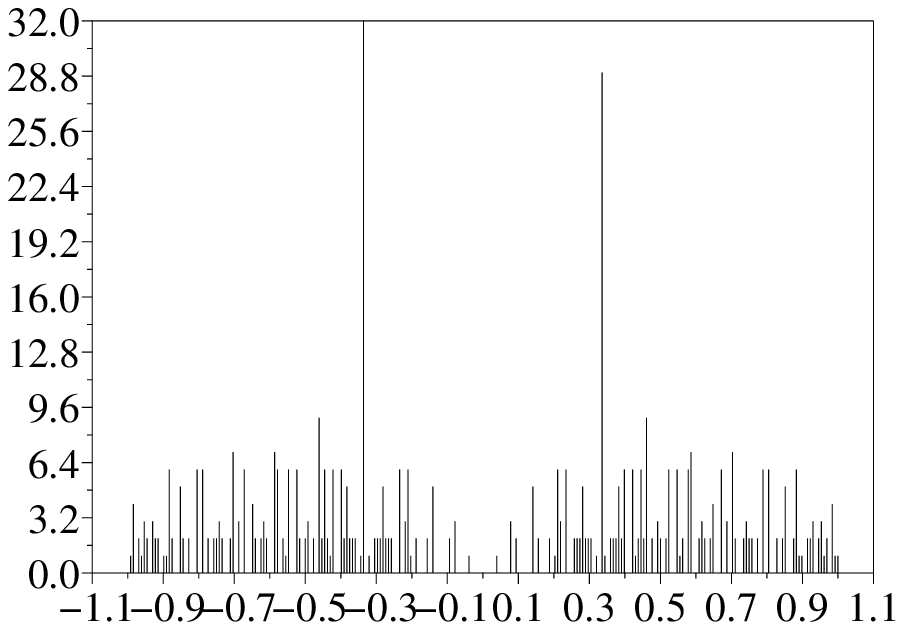,height=150pt}\\
\end{tabular}
\vspace{.3cm}

\noindent\begin{tabular}{@{}p{172pt}p{174pt}@{}} \textbf{Automaton
number $2369$} \vspace{.2cm}

\begin{tabular}{@{}p{48pt}p{200pt}}

$a=\sigma(b,a)$

$b=\sigma(c,a)$

$c=(c,a)$& Group:

Contracting: \emph{no}

Self-replicating: \emph{yes}\\
\end{tabular}

\begin{minipage}{230pt}
\setlength{\fboxrule}{0cm}
\begin{window}[5,r,{\fbox{\shortstack{\hspace{1.6cm}~\phantom{a}\\ \vspace{3.8cm}\\ \phantom{a}}}},{}]
Rels: $a^{-1}ba^{-1}b$, $b^{-1}cb^{-1}c$, $b^{-1}ca^{-1}ba^{-1}c$,\\
$a^{-2}b^{2}a^{-1}b^{-1}ab$, $a^{-2}bab^{-2}ab$,
$a^{-1}ca^{-1}ca^{-1}ca^{-1}c$,\\ $b^{-1}ca^{-1}cb^{-1}ca^{-1}c$,
$a^{-1}bab^{-1}a^{-1}b^{2}a^{-1}$,\\
$a^{-2}bac^{-1}bc^{-1}b^{-1}ab$,
$a^{-1}bc^{-1}ac^{-1}ac^{-1}ac^{-1}b$,\\
$b^{-1}ca^{-2}b^{2}a^{-1}b^{-1}ac$, $b^{-1}ca^{-2}bab^{-2}ac$,\\
$b^{-1}cab^{-2}aba^{-2}c$, $b^{-1}cab^{-1}a^{-1}b^{2}a^{-2}c$,\\
$a^{-1}bab^{-1}c^{-1}bc^{-1}aba^{-1}$
\\
SF: $2^0$,$2^{1}$,$2^{3}$,$2^{7}$,$2^{13}$,$2^{25}$,$2^{47}$,$2^{90}$,$2^{176}$\\
Gr: 1,7,33,143,602,2514\\
\end{window}
\end{minipage}
& \hfill~

\hfill
\begin{picture}(1450,1090)(0,130)
\put(200,200){\circle{200}} 
\put(1200,200){\circle{200}}
\put(700,1070){\circle{200}}
\allinethickness{0.2mm} \put(45,280){$a$} \put(1280,280){$b$}
\put(820,1035){$c$}
\put(164,165){$\sigma$}  
\put(1164,165){$\sigma$}  
\put(664,1022){$1$}       
\put(100,100){\arc{200}{0}{4.71}}     
\path(46,216)(100,200)(55,167)        
\put(300,200){\line(1,0){800}} 
\path(1050,225)(1100,200)(1050,175)   
\spline(287,150)(700,0)(1113,150)     
\path(325,109)(287,150)(343,156)      
\spline(750,983)(1150,287)     
\path(753,927)(750,983)(797,952)      
\spline(650,983)(250,287)      
\path(297,318)(250,287)(253,343)      
\put(700,1211){\arc{200}{2.36}{0.79}} 
\path(820,1168)(771,1141)(779,1196)   
\put(680,240){$_0$} 
\put(193,10){$_1$}  
\put(890,585){$_0$} 
\put(680,77){$_1$}   
\put(545,1261){$_0$}  
\put(460,585){$_1$}  
\end{picture}

\vspace{.3cm}

\hfill\epsfig{file=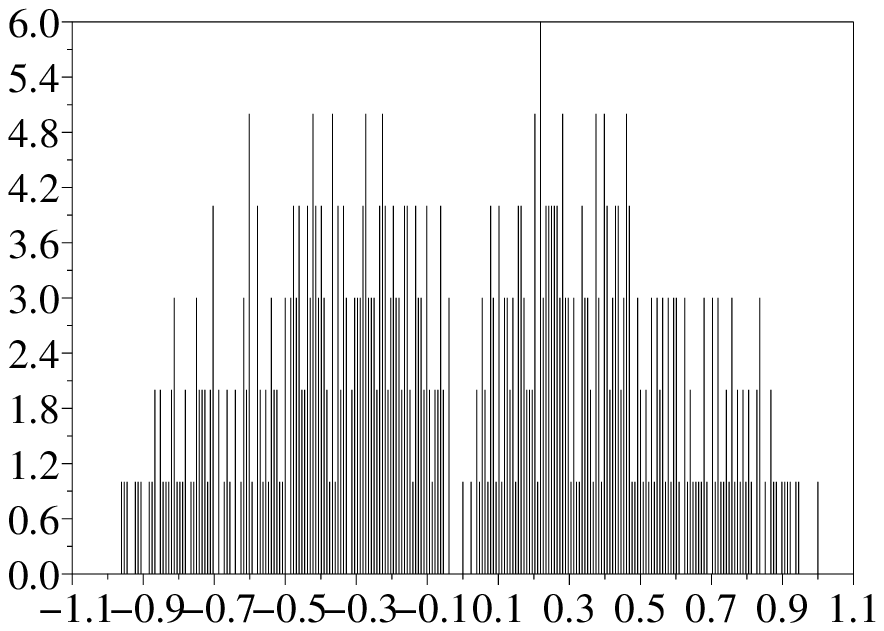,height=150pt}\\
\hline
\end{tabular}
\vspace{.3cm}

\noindent\begin{tabular}{@{}p{172pt}p{174pt}@{}} \textbf{Automaton
number $2851$} \vspace{.2cm}

\begin{tabular}{@{}p{48pt}p{200pt}}

$a=\sigma(a,c)$

$b=\sigma(b,a)$

$c=(c,c)$& Group: Isomorphic to $G_{929}$

Contracting: \emph{no}

Self-replicating: \emph{yes}\\
\end{tabular}

\begin{minipage}{230pt}
\setlength{\fboxrule}{0cm}
\begin{window}[5,r,{\fbox{\shortstack{\hspace{1.6cm}~\phantom{a}\\ \vspace{3.8cm}\\ \phantom{a}}}},{}]
Rels: $c$, $a^{-4}bab^{-1}a^{2}b^{-1}ab$
\\
SF: $2^0$,$2^{1}$,$2^{3}$,$2^{6}$,$2^{12}$,$2^{23}$,$2^{45}$,$2^{88}$,$2^{174}$,$2^{345}$\\
Gr: 1,5,17,53,161,485,1445\\
\end{window}
\end{minipage}
& \hfill~

\hfill
\begin{picture}(1450,1090)(0,130)
\put(200,200){\circle{200}} 
\put(1200,200){\circle{200}}
\put(700,1070){\circle{200}}
\allinethickness{0.2mm} \put(45,280){$a$} \put(1280,280){$b$}
\put(820,1035){$c$}
\put(164,165){$\sigma$}  
\put(1164,165){$\sigma$}  
\put(664,1022){$1$}       
\put(100,100){\arc{200}{0}{4.71}}     
\path(46,216)(100,200)(55,167)        
\spline(200,300)(277,733)(613,1020)   
\path(559,1007)(613,1020)(591,969)    
\spline(287,150)(700,0)(1113,150)     
\path(325,109)(287,150)(343,156)      
\put(1300,100){\arc{200}{4.71}{3.14}} 
\path(1345,167)(1300,200)(1354,216)     
\put(700,1211){\arc{200}{2.36}{0.79}} 
\path(820,1168)(771,1141)(779,1196)   
\put(190,10){$_0$}  
\put(230,700){$_1$} 
\put(1155,10){$_0$}  
\put(680,77){$_1$}   
\put(465,1261){$_{0,1}$}  
\end{picture}

\vspace{.3cm}

\hfill\epsfig{file=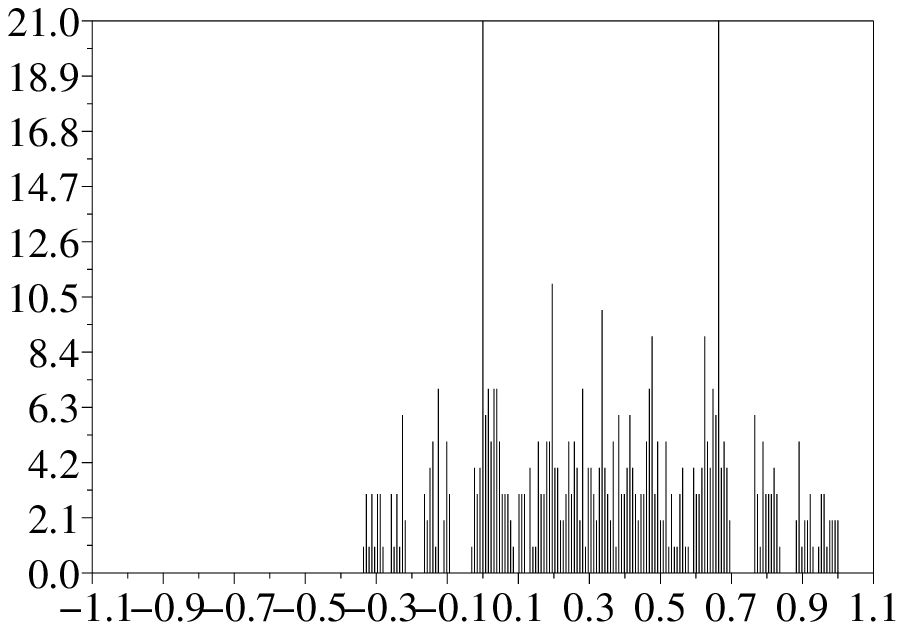,height=150pt}\\
\end{tabular}
\vspace{.3cm}

\noindent\begin{tabular}{@{}p{172pt}p{174pt}@{}} \textbf{Automaton
number $2853$} \vspace{.2cm}

\begin{tabular}{@{}p{48pt}p{200pt}}

$a=\sigma(c,c)$

$b=\sigma(b,a)$

$c=(c,c)$& Group: $IMG\left(\bigl(\frac{z-1}{z+1}\bigr)^2\right)$

Contracting: \emph{yes}

Self-replicating: \emph{yes}\\
\end{tabular}

\begin{minipage}{230pt}
\setlength{\fboxrule}{0cm}
\begin{window}[5,r,{\fbox{\shortstack{\hspace{1.6cm}~\phantom{a}\\ \vspace{3.8cm}\\ \phantom{a}}}},{}]
Rels: $c$, $a^{2}$, $ab^{-1}ab^{-2}ab^{-1}abab^{2}ab$
\\
SF: $2^0$,$2^{1}$,$2^{2}$,$2^{3}$,$2^{5}$,$2^{8}$,$2^{14}$,$2^{25}$,$2^{47}$\\
Gr: 1,4,10,22,46,94,190,375,731,1422,2752,5246,9908\\
Limit space:
\end{window}
\vspace{.1cm} \noindent\epsfig{file=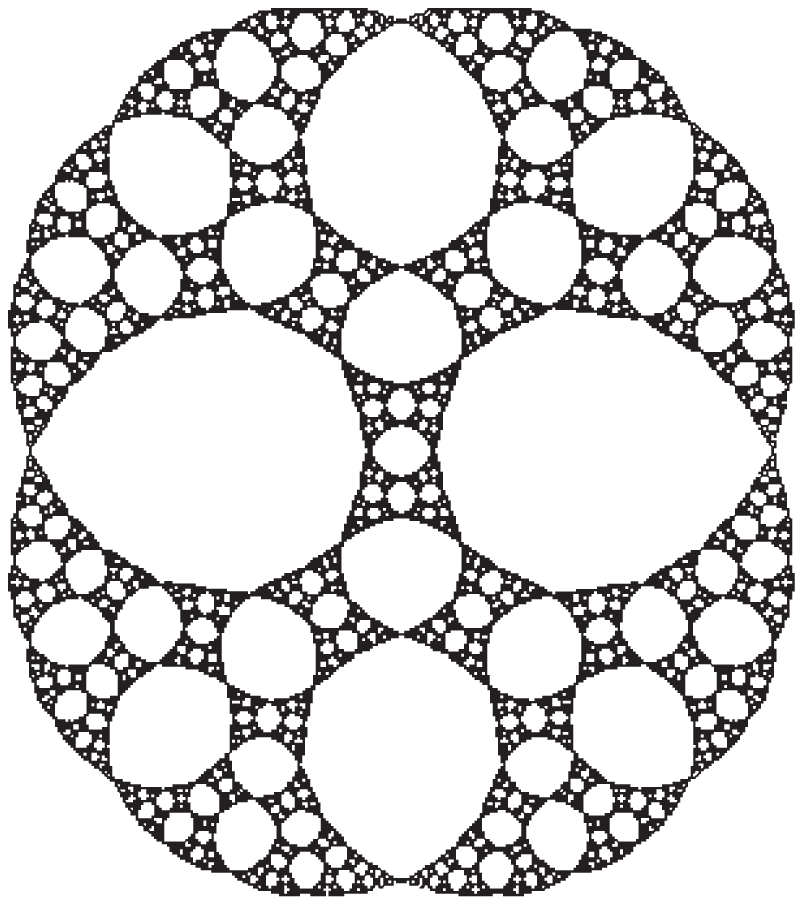,height=135pt}
\end{minipage}
& \hfill~

\hfill
\begin{picture}(1450,1090)(0,130)
\put(200,200){\circle{200}} 
\put(1200,200){\circle{200}}
\put(700,1070){\circle{200}}
\allinethickness{0.2mm} \put(45,280){$a$} \put(1280,280){$b$}
\put(820,1035){$c$}
\put(164,165){$\sigma$}  
\put(1164,165){$\sigma$}  
\put(664,1022){$1$}       
\spline(200,300)(277,733)(613,1020)   
\path(559,1007)(613,1020)(591,969)    
\spline(287,150)(700,0)(1113,150)     
\path(325,109)(287,150)(343,156)      
\put(1300,100){\arc{200}{4.71}{3.14}} 
\path(1345,167)(1300,200)(1354,216)     
\put(700,1211){\arc{200}{2.36}{0.79}} 
\path(820,1168)(771,1141)(779,1196)   
\put(150,700){$_{0,1}$} 
\put(1155,10){$_0$}  
\put(680,77){$_1$}   
\put(465,1261){$_{0,1}$}  
\end{picture}

\vspace{.3cm}

\hfill\epsfig{file=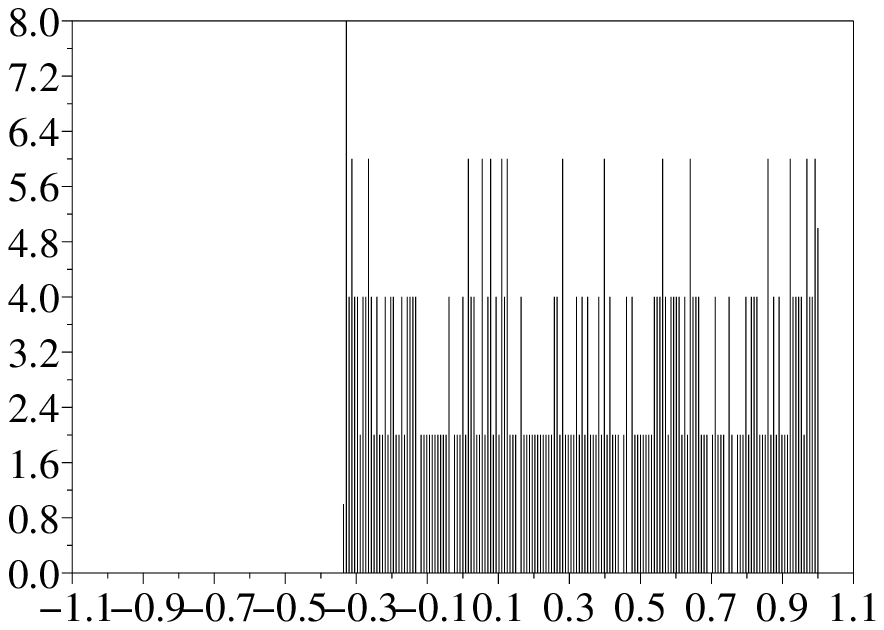,height=150pt}\\
\hline
\end{tabular}
\vspace{.3cm}

\section{Proofs of some facts about the selected groups}
\label{sec_proofs}

We start this section with a few useful observations, which simplify
computations and were used often in the classification process.

First, we need to mention the classification of the groups generated
by $2$-state automata over a $2$-letter alphabet. The following
theorem is proved in~\cite{gns00:automata}.

\begin{theorem}
\label{thm_class2states} There are, up to isomorphism, $6$ different
groups generated by $2$-state automata over a $2$-letter alphabet
(automata of complexity $(2,2)$). Namely, trivial group, $C_2$,
$C_2\times C_2$, infinite cyclic group $\Z$, infinite dihedral group
$D_\infty$ and the lamplighter group $\Z\wr C_2$.
\end{theorem}

The following proposition allows sometimes to see directly from the
automaton structure that the corresponding group is not a torsion
group.

\begin{prop}\label{nontors}
Let $G$ be a group generated by an automaton $A$ over the $2$-letter
alphabet $X=\{0,1\}$ that has the following property. The set of
states of $A$ splits into two nonempty parts $P$ and $Q$ such that
\begin{enumerate}
\item[$(i)$]
one of the parts contains all active states and the other contains
all inactive states;
\item[$(ii)$]
for each state in part $P$, both arrows go to states in the same
part (either both to $P$ pr both to $Q$);
\item[$(iii)$]
for each state in part $Q$, one arrow goes to part $P$ and the other
to part $Q$
\end{enumerate}
Then any element of the group that can be written as a product of
odd number of active generators and odd number of inactive
generators in any order, has infinite order. In particular, $G$ is
not a torsion group.
\end{prop}

\begin{proof}
Let $g$ be such an element. Let us prove by induction on $n$ that,
for each $n \geq 0$, there exists a vertex $v$ fixed by $g^{2^n}$
such that the section of $g^{2^n}$ at $v$ has the same form (i.e. is
a product of odd number of active generators and odd number of
inactive generators).

For $n=0$ this is true. Suppose it is true for $n=k$ and let $h$ be
a section of $g^{n^k}$ at some vertex $v\in X^*$ satisfying the
conditions of the assumption. Since $h$ is a product of odd number
of active states and odd number of inactive states we can write
$h=\sigma(h_1,h_2)$. Since $v$ is fixed under $g^{2^k}$ we have that
$v0$ is fixed under $(g^{2^k})^2=g^{2^{k+1}}$ and
$g^{2^{k+1}}\bigl|_{v0}=h_2h_1$. The element $h_2h_1$ is product (in
some order) of the first level sections of the generators (and/or
their inverses) used to express $h$. By assumption, among these
generators, there are odd number of active and odd number of
inactive. The generators from part $P$, by condition ($ii$), will
produce even number of active and even number of inactive
generators, while the generators from part $Q$, by condition
$(iii)$, will give odd number of generators from both categories,
which proves the induction step. Thus $g^{2^n}\neq 1$ for all $n$.
\end{proof}

There is an algorithm which determines whether a given element of a
self-similar group generated by a finite automaton over the
$2$-letter alphabet $X=\{0,1\}$ acts level transitively on the tree.

The abelianization of $\Aut X^*$ is isomorphic to the infinite
Cartesian product $\prod_{i=0}^\infty C_2$. The canonical
isomorphism sends $g\in G$ to $(c_i \mod 2)_{i=0}^\infty$, where
$c_i$ is the number of vertices $v\in X^i$, such that $g_v$ acts
nontrivially on the first level (i.e.~$c_i$ is the number of active
sections of $g$ at level $i$)

The abelianization group $\prod_{i=0}^\infty C_2$ can be endowed
with the structure of a ring of formal power series $C_2[[t]]$ by
$(a_i)_{i=0}^\infty \mapsto \sum_{i=0}^\infty a_it^i$, where $a_i\in
C_2$.

\begin{algorithm}[Element transitivity]
\label{alg_trans} Let $G$ be a self-similar group generated by a
finite automaton over the $2$-letter alphabet $X=\{0,1\}$, and let
$g\in G$ be given as a product of states and their inverses. Denote
by $\bar g$ the image of $g$ in $C_2[[t]]$. The element $g$ acts
level transitively on $X^*$ if and only if $\bar g=(1,1,1,\ldots)$.

Suppose $g=\sigma^i(g_0,g_1)$, $i\in 0,1$. Then
\[\bar g=i+t\cdot(\overline{g_0} + \overline{g_1}).\]
We can produce similar equations for the sections $g_0$, $g_1$ and
so on. Since $G$ is generated by a finite automaton the number of
different sections of $g$ is finite. Therefore we get finite linear
system of equations over the variables $\{g_v, v\in X^*\}$, whose
solution will express $\bar g$ as a rational function $P(t)/Q(t)$,
where $P,Q$ are polynomials with degrees not higher than $k=|\{g_v,
v\in X^*\}|$.

Expanding this rational function as a power series will produce a
preperiodic sequence of coefficients from $C_2$ with period and
preperiod no longer than $2^k$. In particular, $g$ acts level
transitively if and only if all $c_i$, $i=1,\ldots,2^{k+1}$ are
equal to $1$.
\end{algorithm}

We often need to show that a given group of tree automorphisms is
level transitive. Here is a very convenient necessary and sufficient
condition for this in the case of a binary tree.

\begin{prop}[Group transitivity]\label{prop:transitivity}
A self-similar group of binary tree automorphisms is level
transitive if and only if it is infinite.
\end{prop}

\begin{proof}
Let $G$ be an infinite self-similar group acting on a binary tree.

Level transitivity clearly implies that $G$ is infinite.

For the converse, let us first prove that all level stabilizers
$\Stg(n)$ are different. For this it suffices to show that for every
$n\geq1$ $\Stg(n-1)\setminus\Stg(n)\neq\varnothing$. Since all
stabilizers have finite indices in $G$ and $G$ is infinite we get
that all of them are infinite.

Let $g\in\Stg(n-1)$ be an arbitrary nontrivial element and
$v=x_1\ldots x_k$ be a word of shortest length (one of them) such
that $g(v)\neq v$ (in other words, $g\in\Stg(k-1)$ and such $k$ is
maximal). Clearly $k\geq n$ and we can consider the section
$h=g_{x_1x_2\ldots x_{k-n}}$ which is an element of $G$ because of
self-similarity. The fact that $g\in\Stg({k-1})$ implies
$h\in\Stg(n-1)$. On other hand $x_1\ldots x_{k-n}x_{k-n+1}\ldots
x_k=v\neq g(v)=g(x_1\ldots x_{k-n})h(x_{k-n+1}\ldots x_k)=x_1\ldots
x_{k-n}h(x_{k-n+1}\ldots x_k)$. Therefore $h(x_{k-n+1}\ldots
x_k)\neq x_{k-n+1}\ldots x_k$, thus $h\notin\Stg(n)$ and we found
the desired element.

Now let us prove transitivity by induction on the level. The section
of any nontrivial element at the vertex where it acts nontrivially
gives transitivity on the first level.

Suppose $G$ acts transitively on level $n$. Let
$h\in\Stg(n)\setminus\Stg(n+1)$ be an arbitrary element and let
$w=x_1\ldots x_n\in X^n$ be one of the words such that
$h(wx)=w\overline{x}$, where $\overline{x}=1-x$. For $u=y_1\ldots
y_{n+1}\in X^{n+1}$, let us find an element $g\in G$ such that
$g(w0)=u$. This will prove the induction step.

By inductive assumption there exists $f\in G$ such that
$f(w)=y_1\ldots y_n$. Suppose $f(w0)=y_1\ldots y_n\tilde{y}_{n+1}$.
Then, if $\tilde y_{n+1}=y_{n+1}$ we are done, otherwise $\tilde
y_{n+1}=\overline{y_{n+1}}$ and, for $g=fh$, we obtain
$g(w0)=f(h(w0))=f(w1)=y_1\ldots y_n
\overline{\overline{y_{n+1}}}=u$.
\end{proof}

Note, that the last proof works also for self-similar subgroups of
the infinitely iterated permutational wreath product
$\wr_{i\geq1}C_d$ (the subgroup of $\Aut(\tree)$ consisting of those
automorphisms of the $d$-ary tree for which the activity at every
vertex is a power of some fixed cycle of length $d$. Also, certain
generalizations of this method could be used in more complex
situations (see, for example, the construction of $G_{2240}\cong
C_2\ast C_2\ast C_2$ in~\cite{nekrash:self-similar}, or proof of
transitivity of Sushchansky groups on a subtree
in~\cite{bondarenko_s:sushch}).

We provide below additional information and proofs about the groups
listed in Section~\ref{sec_tables}. \smallskip

\noindent\textbf{739}: $C_2\ltimes(C_2\wr\Z)$. For $G_{739}$, we
have $a=\sigma(a,a)$, $b=(b,a)$ and $c=(a,a)$.

All generators have order $2$. The elements $u=acba=(1,ba)$ and
$v=bc=(ba,1)$ generate $\Z^2$ because $ba=\sigma(1,ba)$ is the
adding machine and has infinite order. Also we have $ac=\sigma$ and
$\langle u,v\rangle$ is normal in $H=\langle u,v,\sigma\rangle$
because $u^\sigma=v$ and $v^\sigma=u$. In other words, $H\cong
C_2\ltimes(\Z\times\Z) =C_2\wr\Z$.

Furthermore, $G_{739}=\langle H,a\rangle$ and $H$ is normal in
$G_{739}$ because $u^a=v^{-1}$, $v^a=u^{-1}$ and $\sigma^a=\sigma$.
Thus $G_{739}=C_2\ltimes(C_2\wr\Z)$, where the action of $C_2$ on
$H$ is specified above.

\noindent\textbf{744}: For $G_{744}$, we have $a=\sigma(c,b)$,
$b=(b,a)$ and $c=(a,a)$.

Since $(a^{-1}c)^2=(c^{-1}ab^{-1}a,b^{-1}ac^{-1}a)$ and
$c^{-1}ab^{-1}a = ((c^{-1}ab^{-1}a)^{-1},a^{-1}c)$ we see that
$(a^{-1}c)^2$ fixes the vertex $01$ and its section at this vertex
is equal to $a^{-1}c$. Hence, $a^{-1}c$ has infinite order.

Furthermore, the element $c^{-1}ab^{-1}a$ has infinite order, fixes
the vertex $00$ and its section at this vertex is equal to
$c^{-1}ab^{-1}a$. Therefore $G_{744}$ is not contracting (all powers
of $c^{-1}ab^{-1}a$ would have to belong to the nucleus).

\noindent\textbf{748}: $D_4\times C_2$. For $G_{748}$, we have
$a=\sigma(a,a)$, $b=(c,a)$ and $c=(a,a)$.

It follows from the relations
$a^{2}=b^{2}=c^{2}=acac=bcbc=abababab=1$ that $G_{748}$ is a
homomorphic image of $D_4\times C_2$. Since $a\neq 1$, $b \neq 1$
and $(ab)^2 \neq 1$, it follows that $\langle a,b \rangle = D_4$.
One can verify directly that $c$ is not equal to any of the four
elements in $\langle a,b \rangle$ that stabilize level 1 (namely
$1$, $b$, $aba$ and $abab$).  Thus $G_{748}=D_4 \times C_2$.

\noindent\textbf{753}: For $G_{753}$, we have $a=\sigma(c,b)$,
$b=(c,a)$ and $c=(a,a)$.

Since $ab^{-1}=\sigma(1,ba^{-1})$, this element is conjugate to the
adding machine.

For a word $w$ in $w\in\{a^{\pm1},b^{\pm1},c^{\pm1}\}^*$, let
$|w|_a$, $|w|_b$ and $|w|_c$ denote the sum of the exponents of $a$,
$b$ and $c$ in $w$. Let $w$ represents the element $g \in G$. If
$|w|_a$ and $|w|_b$ are odd, then $g$ acts transitively on the first
level, and $g^2|_0$ is represented by a word $w_0$, which is the
product (in some order) of all first level sections of all
generators appearing in $w$. Hence, $|w_0|_a=|w|_b+2|w|_c$ and
$|w_0|_b=|w|_a$ are odd again. Therefore, similarly to
Proposition~\ref{nontors}, any such element has infinite order.

In particular $c^2ba$ has infinite order. Since
$a^4=(caca,a^4,acac,a^4)$ and $caca=(baca,c^2ba,bac^2,caba)$, the
element $a^4$ has infinite order (and so does $a$). Since $a^4$
fixes the vertex $01$ and its section at that vertex is equal to
$a^4$, the group $G_{753}$ is not contracting.

\noindent\textbf{771}: $\mathbb Z^2$. For $G_{771}$, we have
$a=\sigma(c,b)$, $b=(b,b)$ and $c=(a,a)$.

Since $G_{771}$ is finitely generated, abelian, and self-replicating
(easy to check), it follows from~\cite{nekrash_s:12endomorph} that
it is free abelian. There are two options: either it has rank $1$ or
rank $2$ (since $b=1$). Let us prove that the rank is $2$. For this
it is sufficient to show that $c^n \neq a^m$ in $G$. Assume on the
contrary that $c^n=a^m$ for some integer $n$ and $m$ and choose such
integers with minimal $|n|+|m|$. Since $c^n$ stabilizes level 1, $m$
must be even and we have $(a^n,a^n)=c^n = a^m = (c^{m/2},c^{m/2})$.
But then $a^n=c^{m/2}$ and by the minimality assumption $m$ must be
$0$, implying $c^n=1$. The last equality can only be true for $n=0$
since $G_{771}$ is torsion free (free abelian) and $c \neq 1$. Thus
$G_{771}\cong\mathbb Z^2$.

\noindent\textbf{775}: $C_2\ltimes
IMG\left(\bigl(\frac{z-1}{z+1}\bigr)^2\right)$. For $G_{775}$, we
have $a=\sigma(a,a)$, $b=(c,b)$, $c=(a,a)$.

We have $a^2=b^2=c^2=1$, $ac=ca=\sigma(1,1)$ and $ba=\sigma(ba,ca)$.
Hence, for the subgroup $H=\langle ba,ca\rangle \leq G$, we have
$H\cong G_{2853}\cong
IMG\left(\bigl(\frac{z-1}{z+1}\bigr)^2\right)$. On the other hand
$H$ is normal in $G$ since $(ba)^a=ab=(ba)^{-1}$ and $(ca)^a=ac=ca$.
Thus $G\cong C_2\ltimes H$, where $C_2$ is generated by $a$ and the
action of $a$ on $H$ is given above. It is proved below (see
$G_{783}$), that $G_{775}\cong G_{783}$. Therefore $G_{775}$ also
contains a torsion free subgroup of index $4$.

\noindent\textbf{783$\cong G_{775}$}: $C_2\ltimes
IMG\left(\bigl(\frac{z-1}{z+1}\bigr)^2\right)$. For $G_{783}$, we
have $a=\sigma(c,c)$, $b=(c,b)$ and $c=(a,a)$.

All generators have order $2$ and $a$ commutes with $c$. Conjugating
this group by the automorphism $\gamma=(c\gamma,\gamma)$ yields an
isomorphic group generated by the 4-state automaton defined by the
recursive relations $a'=\sigma(1,1)$, $b'=(c',b')$ and $c'=(a',a')$.
On the other hand, we obtain the same automaton after conjugating
$G_{775}$ by $\mu=(a\mu,\mu)$ (here $a$ denotes the generator of
$G_{775}$).

It can be proved that the subgroup $H=\langle ba,cabc\rangle$ is
torsion free and not metabelian. Furthermore, $G_{783}=\langle
a,c\rangle\ltimes H\cong (C_2\times C_2)\ltimes H$. The group
$G_{783}$ is regular weakly branch group over $H''$.

Since $bca=\sigma(bca,a)$, $G_{783}=\langle acb,a,c\rangle\cong
G_{2205}$.

\noindent\textbf{803$\cong G_{771}$}: $\mathbb Z^2$. For $G_{803}$,
we have $a=\sigma(b,a)$, $b=(c,c)$, $c=(a,a)$.

Since $G_{771}$ is finitely generated, abelian, and
self-replicating, it follows from~\cite{nekrash_s:12endomorph} that
it is free abelian. Consider the $\frac12$-endomorphism
$\phi:\mathop{\rm Stab}\nolimits_{G_{803}}(1)\to G_{803}$ associated
to the vetrex $0$, given by $\phi(g)=h$ for $g\in\mathop{\rm
Stab}\nolimits_{G_{803}}(1)$, provided $g=(h,*)$. Consider also the
linear map $A:\mathbb C^3\to\mathbb C^3$ induced by $\phi$. It has
the following matrix representation with respect to the basis
corresponding to the triple $\{a,b,c\}$:
\[
 A=\left(
 \begin{array}{ccc}
  \frac12 & 0 & 1 \\
  \frac12 & 0 & 0 \\
  0 & 1 & 0
\end{array}
\right).
\]

Its characteristic polynomial
$\chi_A(\lambda)=-\lambda^3+\frac12\lambda^2+\frac12$ has three
distinct complex roots $\lambda_1=1$,
$\lambda_2=-\frac14-\frac14i\sqrt7$ and
$\lambda_3=-\frac14+\frac14i\sqrt7$. Choose an eigenvector $v_i$
associated to the eigenvalue $\lambda_i$, $i=1,2,3$. In particular,
we may choose $v_1=(2,1,1)$, which shows that $a^2bc=1$ in
$G_{803}$. In order to show that $a^{2m}c^n\neq1$ (except when
$m=n=0$) we will prove that the vector $v=(2m,0,n)$ is eventually
pushed out from the domain corresponding to the first level
stabilizer, i.e. from the set $D=\{(2i,j,k), i,j,k\in\mathbb Z\}$,
by iterations of the action of $A$.

Consider the expansion of $v$ in the basis $\{v_1,v_2,v_3\}$:
$v=a_1v_1+a_2v_2+a_3v_3$. Since $m \neq 0$ or $n \neq 0$, $v$ is not
a scalar multiple of $v_1$. We have
$A^t(v)=a_1v_1+\lambda_2^ta_2v_2+\lambda_3^ta_3v_3\to a_1v_1$, as
$t\to\infty$, since $|\lambda_2|=|\lambda_3|<1$. We can choose a
neighborhood of $a_1v_1$ that does not contain points from $D$,
except maybe $a_1v_1$. Eventually $A^t(v)$ will be in this
neighborhood and, since $A^t(v)\neq a_1v_1$ for all $t$, $A^t(v)$
will be outside of $D$. This implies that the word $a^{2m}c^n$
represents a nontrivial element in $G_{803}$. Thus $G_{803} =
\langle a,c \rangle \cong \mathbb Z^2$.

\noindent\textbf{846}: $C_2\ast C_2\ast C_2$. This is a result of
Muntyan. See the proof in~\cite{nekrash:self-similar}. In
particular, $G_{846}$ contains a self-similar free group of rank $2$
generated by a 6-state automaton. The automaton 846 is sometimes
called Bellaterra automaton.

\noindent\textbf{852}: Basilica group $\B = IMG(z^2-1)$. First
studied in~\cite{grigorch_z:basilica}, where it is proved that $\B$
is not in the class $SG$ of sub-exponentially amenable groups, does
not contain a free subgroup of rank $2$, and the nontrivial
generators $a$ and $b$ generate a free subsemigroup. Spectral
properties are considered in~\cite{grigorch_z:basilica_sp}. It is
proved in~\cite{bartholdi-v:basilica} that $\B$ is amenable,
providing the first example of an amenable group not in the class
$SG$.

\noindent\textbf{857}: For $G_{857}$, we have $a=\sigma(b,a)$,
$b=(c,b)$ and $c=(b,a)$.

Let us prove that $b$ has infinite order. For any $w\in X^*$
$b(w0^\infty)=b(w0)b_{w0}(0^\infty)$. Since $b_{w0}$ equals either
$b$ or $c$ and $b(0^\infty)=c(0^\infty)=0^\infty$, we have
$b(w0^\infty)=b(w0)0^\infty$. Therefore all elements in the forward
orbit of $010^\infty$ under the action of $b$ end in $0^\infty$. The
length of the non-zero prefix of any infinite word ending in
$0^\infty$ cannot decrease under the action of $b$. Indeed, for any
$w\in X^*$ $b(w10^\infty)=b(w)b_{w}(10^\infty)$. The section $b_w$
is one of the three generators, for which we have
$a(10^\infty)=010^\infty$, $b(10^\infty)=10^\infty$ and
$c(10^\infty)=110^\infty$.

On the other hand, the length of the non-zero prefix along the orbit
cannot stabilize, because in this case the orbit must be finite and
we must have $b^k(010^\infty)=010^\infty$, for some $k\geq1$. But
this is impossible since $b(010^\infty)=0110^\infty$ and thus the
length of the non-zero prefix of $b^k(010^\infty)$ must be at least
$3$. Thus the orbit is infinite and $b$ has infinite order.

Since $b=(c,b)$, $G_{857}$ is not contracting.

\noindent\textbf{858}: For $G_{858}$, we have $a=\sigma(c,a)$,
$b=(c,b)$ and $c=(b,a)$.

The element $ab^{-1}=\sigma(1,ab^{-1})$ is the adding machine.

Using the same approach as for $G_{857}$ one can show that $c$ has
infinite order. Namely the length of the non-zero prefix of the
forward orbit of $10^\infty$ under $c$ is nondecreasing, which then
implies that this orbit is not finite.

Since $b=(c,b)$, $G_{858}$ is not contracting.

\noindent\textbf{870}: Baumslag-Solitar group $BS(1,3)$. For
$G_{870}$, we have  $a=\sigma(c,b)$, $b=(a,c)$, and $c=(b,a)$.

The automaton satisfies the conditions of Proposition~\ref{nontors}.
Thus, $ab$ has infinite order in $G$, which implies that
$bc=(ab,ca)$, $a^2=(bc,cb)$ also have infinite order. Hence, we can
claim the same for $a$ and $b=(a,c)$.

Furthermore, the element
$\mu=b^{-1}a=\sigma(1,a^{-1}b)=\sigma(1,\mu^{-1})$ also has infinite
order (it is conjugate of the adding machine). Since
$a^{-1}c=\sigma(1,c^{-1}a)=\mu$ we have $c=ab^{-1}a$ and $G=\langle
a,b\rangle=\langle \mu,b\rangle$. Let us check that $b^{-1}\mu
b=\mu^3$. Since $b^{-1}\mu b = \sigma(c^{-1}a,a^{-2}bc)$ and $\mu^3
= \sigma(\mu^{-1},\mu^{-2})$ all we need to check is that
$a^{-2}bc=a^{-1}ba^{-1}b$, i.e.~$a^{-1}bcb^{-1}ab^{-1}=1$. The last
is correct since $a^{-1}bcb^{-1}ab^{-1}=(1,b^{-1}aba^{-1}bc^{-1})$
and $b^{-1}aba^{-1}bc^{-1}$ is a conjugate of the inverse of
$a^{-1}bcb^{-1}ab^{-1}$. Thus, since $b$ and $\mu$ have infinite
order, $G_{870}\cong BS(1,3)$.

See~\cite{bartholdi_s:bsolitar} for realizations of $BS(1,3)$ and
other Baumslag-Solitar groups by automata.

\noindent\textbf{878}: $C_2\ltimes IMG(1-\frac1{z^2})$. For
$G_{878}$, we have $a=\sigma(b,b)$, $b=(b,c)$ and $c=(b,a)$.

Denote $x=bc$ and $y=ca$. All generators have order $2$, and
therefore the subgroup $H=\langle x,y\rangle$ is a normal subgroup
of index $2$ in $G_{878}$. Moreover $G_{878}\cong C_2\ltimes H$,
where $C_2$ is generated by $c$ and the action of $C_2$ on $H$ is
given by $x^c = x^{-1}$ and $y^c = y^{-1}$. We have $x=(1,ca)=(1,y)$
and $y=\sigma(ab,1)=\sigma(y^{-1}x^{-1},1)$. Exchanging the letters
$0$ and $1$ leads to an isomorphic copy of $H$ defined by $x=(y,1)$
and $y=\sigma(1,y^{-1}x^{-1})$, which is the iterated monodromy
group $IMG(1-\frac1{z^2})$, according to~\cite{bartholdi_n:rabbit}.
Thus, $G_{878}\cong C_2\ltimes IMG(1-\frac1{z^2})$.

\noindent\textbf{929}. See $G_{2851}$.

\noindent\textbf{942}: For $G_{942}$, we have $a=\sigma(c,b)$,
$b=(c,b)$ and $c=(c,a)$.

It is known~\cite{grigorchuk-z:l2} that the group $L=\langle
a',b'\rangle$ defined by
\begin{alignat*}{2}
  &a' = \sigma &&(a',b')\\
  &b' =        &&(a',b')
\end{alignat*}
is the lamplighter group $\Z \wr C_2$ (compare to~\eqref{lampl}).
Consider the subtree $Y^*$ of $X^*$ consisting of all words over the
alphabet $Y=\{01,11\}$. The element $a$ swaps the letters of $Y$ and
$b$ fixes them. Since $a_{01}=b_{01}=a$, $a_{11}=b_{11}=b$, the tree
$Y^*$ is invariant under the action of $H=\langle a,b\rangle$ and
the action of $H$ on $Y^*$ coincides with the action of the
lamplighter group $L=\langle a',b' \rangle$ on $X^*$ (with the
identification $0\leftrightarrow 01, \ 1 \leftrightarrow 11$).
Therefore the map $\phi:H \to L$ given by $a \mapsto a'$, $b \mapsto
b'$ extends to a homomorphism. We claim that this homomorphism has
trivial kernel. Indeed, let $w=w(a,b)$ be a group word representing
an element of the kernel of $\phi$. Since the word $w(a',b')$
represents the identity in $L$ the total exponent of $a$ in $w$ must
be even and the total exponent $\epsilon$ of both $a$ and $b$ in $w$
must be 0. But in that case the element $g=w(a,b)$ fixes the top two
levels of the tree $X^*$ and has decomposition
\[ g = (c^\epsilon,*,c^\epsilon,*), \]
where the $*$'s denote words over $a$ and $b$ representing the
identity in $H$ (these words correspond to the first level sections
of $w(a',b')$ in $L$). Therefore $g=1$ and the kernel of $\phi$ is
trivial.

Thus, the lamplighter group is a subgroup of $G_{942}$.

\noindent\textbf{968}. For $G_{968}$, we have $a=\sigma(b,b)$,
$b=(c,c)$ and $c=(c,a)$.

This group contains $\Z^5$ as a subgroup of index $16$. It is
contracting with nucleus consisting of $73$ elements, whose
self-similar closure consists of $77$ elements. All generators have
order $2$.

Let $x=(ac)^2$, $y=bcba$, and consider the subgroup $K=\langle x,y
\rangle$. Direct computations show that $x$ and $y$ commute
($xy=\sigma(bcbacacb,ba)$ and $yx=\sigma(cacabc,ba)$). Conjugating
by $\gamma=(b\gamma,a\gamma)$ leads to the self-similar copy $K'$ of
$K$ generated by $x' = ((y')^{-1}, (y')^{-1})$ and $y =
\sigma(x',y')$, where $x' = x^\gamma$ and $y'=y^\gamma$. Since
$(y')^2 = (x'y',x'y')$, the virtual endomorphism of $K'$ is given by
the matrix
\[
 A=\left(\begin{array}{cc}
               0&\frac12\\
              -1&\frac12\\
 \end{array} \right).
\]
The eigenvalues $\lambda=\frac14\pm\frac14\sqrt{7}i$ of this matrix
are not algebraic integers, hence, according
to~\cite{nekrash_s:12endomorph}, the group $K'$ is free abelian of
rank $2$, and so is $K$.

Since all generators have order $2$, the subgroup $H=\langle
ab,bc\rangle$ has index $2$ in $G_{968}$. The stabilizer $\St_H(2)$
of the second level has index $8$ in $H$. Moreover, the quotient
group is isomorphic to the dihedral group $D_4$ (since $ab$ acts on
the second level by permuting $00 \leftrightarrow 10$ and $01
\leftrightarrow 11$, while $bc$ acts by permuting $10
\leftrightarrow 11$). The stabilizer $\St_H(2)$, conjugated by the
element $g=(b,c,b,1)$, is generated by
\[
\begin{array}{lllllll}
g_1=\left((bc)^2\right)^g&=(bcbc)^g=&(1,&1,&y,&y^{-1}&),\\
g_2=\left((bc)^2\right)^{bag}&=(acbcba)^g=&(y,&y,&1,&1&),\\
g_3=\left((ab)^2\right)^{bcg}&=(cbabac)^g=&(1,&x,&x,&1&),\\
g_4=\left((ab)^2\right)^g&=(abab)^g=&(1,&x,&1,&x^{-1}&),\\
g_5=\left((ab)^2\right)^{(bc)^{ba}g}&=(abcbabacba)^g=&(x,&1,&1,&x^{-1}&),\\
g_6=\left((ab)^2\right)^{bc(bc)^{ba}g}&=(abcacbabacacba)^g=&(x,&1,&x,&1&).\\
\end{array}
\]
Therefore, all $g_i$ commute and $g_6=g_5g_3g_4^{-1}$. If
$\prod_{i=1}^5g_i^{n_i}=1$, then all sections must be trivial,
hence,
$x^{n_5}y^{n_2}=x^{n_3+n_4}y^{n_2}=x^{n_3}y^{n_1}=x^{n_4+n_5}y^{n_1}=1$.
But $K$ is free abelian, whence $n_i=0$, $i=1,\ldots,5$. Thus,
$\St_H(2)$ is a free abelian group of rank $5$.

\noindent\textbf{2205=$G_{775}$}. $C_2\ltimes
IMG\left(\bigl(\frac{z-1}{z+1}\bigr)^2\right)$. For $G_{2205}$, we
have $a=\sigma(c,c)$, $b=\sigma(b,a)$ and $c=(a,a)$. See $G_{783}$
for an isomorphism.

\noindent\textbf{2212}: Klein bottle group, $\langle a,b\ \bigl|\
a^2=b^2\rangle$, contains $\Z^2$ as subgroup of index 2. For
$G_{2212}$, we have $a=\sigma(a,c)$, $b=\sigma(c,a)$ and $c=(a,a)$.

Since $ac=\sigma(a^2,ca)$ and $ca=\sigma(a^2,ac)$, the generators
$a$ and $c$ commute. Further, $a^2c = (ca^2,aca) = (a^2c,a^2c)$,
which shows that $c=a^{-2}$, and therefore $a=\sigma(a,a^{-2})$ and
$b=\sigma(a^{-2},a)$. Since $a^2=(a^{-1},a^{-1})$, the element $a$
has infinite order and so does $x=ab^{-1}=(a^3,a^{-3})$. Finally,
since $x^a = b^{-1}a = (a^{-3},a^3) = x^{-1}$, we have $G_{2212} =
\langle x,a \mid x^a = x^{-1}\rangle$ and $G_{2212}$ is the Klein
bottle group. Going back to the generating set consisting of $a$ and
$b$, we get the presentation $G_{2212} = \langle a,b \mid a^2 =
b^2\rangle$.

\noindent\textbf{2240}: Free group of rank $3$. The automaton
generating this group first appeared in~\cite{aleshin:free}. It is
proved in~\cite{vorobets:aleshin} that $G_{2240}$ is a free group of
rank 3 with basis $\{a,b,c\}$. This is the smallest example of a
free nonabelian group among all automata over a $2$-letter alphabet
(see Theorem~\ref{thm:class_free}).

\noindent\textbf{2277}: $C_2\ltimes (\Z\times\Z)$. For $G_{2277}$,
we have $a=\sigma(c,c)$, $b=\sigma(a,a)$, $c=(b,a)$.

All generators have order $2$. Let $x=cb$ and $y=ab$ and let
$H=\langle x,y \rangle$. Then $x=\sigma(1,y^{-1})$ and
$y=(xy^{-1},xy^{-1})$. It is easy to check that $x$ and $y$ commute
and that $H$ is self-replicating. The matrix of the associated
virtual endomorphism is given by
\[
 A =
 \begin{pmatrix}
   0 & 1 \\ -1/2 & -1
 \end{pmatrix}.
\]
Since the eigenvalues $-\frac12 \pm \frac12 i$ are not algebraic
integers, according to~\cite{nekrash_s:12endomorph} $H$ is free
abelian of rank 2.

The subgroup $H$ is normal of index 2 in $G_{2277}$ because the
generators of $G_{2277}$ are of order $2$. Thus $G_{2277}=\langle
H,b\rangle = C_2\ltimes(\Z\times\Z)$, where the action of
$C_2=\langle b \rangle$ on $H$ is by inversion of the generators.

\noindent\textbf{2369}. For $G_{2369}$ we have $a=\sigma(b,a)$,
$b=\sigma(c,a)$ and $c=(c,a)$.

For any vertex $v\in X^*$, we have $a_{v1}=a$, $a_{v10}=b$ and
$a_{v10^{n+2}}=c$, for $n\geq 0$. Therefore, for any vertex $w\in
X^*$, $a(w10^\infty)=a(w1)110^\infty$ and the forward orbit of
$10^\infty$ under $a$ is infinite, because the length of the
non-zero prefix grows by 2 with each application of $a$. Thus $a$
has infinite order.

Since $a^2=(ab,ba)$, the element $ab$ also has infinite order.
Furthermore, $ab=(ac,ba)$ and $ba=(ab,ca)$. Thus, $G_{2369}$ is not
contracting.

\noindent\textbf{2851=$G_{929}$}. For $G_{2851}$ we have
$a=\sigma(a,1)$, $b=\sigma(b,a)$, $c=(c,c)=1$.

The element $a$ is conjugate of the adding machine (in fact it is
its inverse). Since $ba^{-1}=(a,ba^{-1})$, the order of $ba^{-1}$ is
infinite and $G_{2851}$ is not contracting.

The group $G_{2851}$ is a regular weakly branch group over $G'$
since it is self-replicating and $[a^2,b]=([a,b],1)$.

The subsemigroup $\langle a,b\rangle$ is free. Indeed, let $w$ be a
nonempty word in $\{a,b\}^*$. If $w=1$ in $G_{2851}$, then $w$
contains both $a$ and $b$, because they both have infinite order.
Suppose the length of $w$ is minimal among all nonempty words over
$\{a,b\}$ representing the identity element in $G_{2851}$. Then one
of the projections of $w$ will be shorter than $w$, nonempty, and
will represent the identity in $G_{2851}$, which contradicts the
minimality assumption. Thus $w\neq 1$ in $G_{2851}$, for any
nonempty word in $\{a,b\}^*$.

Now let $w$ and $v$ be two words in $\{a,b\}^*$ with minimal sum
$|w|+|v|$ such that $w=v$ in $G_{2851}$. Suppose $w$ ends in $a$ and
$v$ ends in $b$. Then
\begin{enumerate}
\item
if $w$ ends in $a^2$ then $w_0$ is a word that is shorter than $w$
ending in $a$, while $v_0$ is a word not longer than $v$ ending in
$b$. Since $w_0=v_0$ in $G_{2851}$ and
$\bigl|w_0\bigr|+\bigl|v_0\bigr|<|w|+|v|$, we have a contradiction.

\item
if $w$ ends in $ba$ then $w_1$ is a word shorter than $w$ ending in
$b$, while $v_1$ is a word not longer than $v$ ending in $a$. Since
$w_1=v_1$ in $G_{2851}$ and
$\bigl|w_1\bigr|+\bigl|v_1\bigr|<|w|+|v|$,  we have a contradiction.

\item
if $w=a$ then $v_1=1$ in $G$ and $v_1$ is a nonempty word, which is
impossible, as already proved above.
\end{enumerate}

Thus $G$ has exponential growth. On the other hand, the orbital
Schreier graph $\Gamma(G,000\ldots)$ has intermediate growth
(see~\cite{benjamini_h:omega_per_graphs,bond_cn:amenable}).

The group $G_{2851}$ coincides with $G_{929}$ as subgroup of $\Aut(
X^*)$. Indeed, $G_{2851}=\langle a^{-1}=\sigma(1,a^{-1}),
b^{-1}a=(b^{-1}a,a^{-1})\rangle=G_{929}$. Therefore all properties
proved for $G_{2851}$ above hold also for $G_{929}$.

\noindent\textbf{2853}:
$IMG\left(\bigl(\frac{z-1}{z+1}\bigr)^2\right)$. For $G_{2853}$, we
have $a=\sigma(c,c)$, $b=\sigma(b,a)$ and $c=(c,c)=1$.

It is proved in~\cite{bartholdi_n:rabbit} that
$IMG\left(\bigl(\frac{z-1}{z+1}\bigr)^2\right)$ is generated by
$\alpha=\sigma(1, \beta)$ and $\beta=(\alpha^{-1}\beta^{-1},
\alpha)$. We have then $\beta\alpha=\sigma(\alpha, \alpha^{-1})$.
Conjugate the right hand side of the wreath recursion by
$(1,\alpha)$ to obtain a copy of
$IMG\left(\bigl(\frac{z-1}{z+1}\bigr)^2\right)$ given by
$\beta=(\alpha^{-1}\beta^{-1}, \alpha)$, $\beta\alpha=\sigma$ and
$\alpha=\sigma(\alpha^{-1}, \beta\alpha)$ (this is equivalent to
conjugating by $\gamma=(\gamma,\alpha\gamma)$ in $\Aut(X^*)$).

This shows that $G_{2853}$ is isomorphic to
$IMG\left(\bigl(\frac{z-1}{z+1}\bigr)^2\right)$ via the isomorphism
$a\mapsto\beta\alpha$ and $b\mapsto\alpha$. Moreover, they are
conjugate by the element $\delta=(\delta_1, \delta_1)$, where
$\delta_1=\sigma(\delta, \delta)$ (this is the automorphism of the
tree changing all letters which stand on even places).

Consequently, the limit space of $G_{2853}$ is the Julia set of the
rational map $z \mapsto \bigl(\frac{z-1}{z+1}\bigr)^2$.

The group $G_{2853}$ is contained in $G_{775}$ as a subgroup of
index $2$ (see $G_{775}$). It contains the torsion free subgroup $H$
mentioned in the discussion of $G_{775}$ as subgroup of index 2 and
is a weakly branch group over $H''$. All Schreier graphs on the
boundary of the tree have polynomial growth of degree $2$. Diameters
of Schreier graphs on the levels grow as $\sqrt{2}^n$
(see~\cite{bond_n:schreier} for details).

\newcommand{\etalchar}[1]{$^{#1}$}
\def\cprime{$'$} \def\cprime{$'$} \def\cprime{$'$} \def\cprime{$'$}
  \def\cprime{$'$}


\begin{thebibliography}{BKNV05}

\bibitem[AHM04]{aitken-h-m:iterated}
Wayne Aitken, Farshid Hajir, and Christian Maire.
\newblock Finitely ramified iterated extensions.
\newblock (preprint), 2004.

\bibitem[Ale72]{aleshin:burnside}
S.~V. Ale{\v{s}}in.
\newblock Finite automata and the {B}urnside problem for periodic groups.
\newblock {\em Mat. Zametki}, 11:319--328, 1972.

\bibitem[Ale83]{aleshin:free}
S.~V. Aleshin.
\newblock A free group of finite automata.
\newblock {\em Vestnik Moskov. Univ. Ser. I Mat. Mekh.}, (4):12--14, 1983.

\bibitem[AV05]{abert-v:dimension}
Mikl{\'o}s Ab{\'e}rt and B{\'a}lint Vir{\'a}g.
\newblock Dimension and randomness in groups acting on rooted trees.
\newblock {\em J. Amer. Math. Soc.}, 18(1):157--192 (electronic), 2005.

\bibitem[BCSN]{bond_cn:amenable}
I.~Bondarenko, T.~Chekerini-Sil{\cprime}bersta{\u\i}n, and
V.~Nekrashevych.
\newblock Amenable graphs with dense holonomy and no compact isometry groups.
\newblock {I}n preparation.

\bibitem[BG00a]{bartholdi_g:spectrum}
L.~Bartholdi and R.~I. Grigorchuk.
\newblock On the spectrum of {H}ecke type operators related to some fractal
  groups.
\newblock {\em Tr. Mat. Inst. Steklova}, 231(Din. Sist., Avtom. i Beskon.
  Gruppy):5--45, 2000.

\bibitem[BG00b]{bartholdi-g:lie}
Laurent Bartholdi and Rostislav~I. Grigorchuk.
\newblock Lie methods in growth of groups and groups of finite width.
\newblock In Michael~Atkinson et~al., editor, {\em Computational and Geometric
  Aspects of Modern Algebra}, volume 275 of {\em London Math. Soc. Lect. Note
  Ser.}, pages 1--27. Cambridge Univ. Press, Cambridge, 2000.

\bibitem[BGN03]{bartholdi-g-n:fractal}
Laurent Bartholdi, Rostislav Grigorchuk, and Volodymyr Nekrashevych.
\newblock From fractal groups to fractal sets.
\newblock In {\em Fractals in Graz 2001}, Trends Math., pages 25--118.
  Birkh\"auser, Basel, 2003.

\bibitem[BG{\v{S}}03]{bartholdi-g-s:branch}
Laurent Bartholdi, Rostislav~I. Grigorchuk, and Zoran
{\v{S}}uni{\'k}.
\newblock Branch groups.
\newblock In {\em Handbook of algebra, Vol. 3}, pages 989--1112. North-Holland,
  Amsterdam, 2003.

\bibitem[BH05]{benjamini_h:omega_per_graphs}
Itai Benjamini and Christopher Hoffman.
\newblock {$\omega$}-periodic graphs.
\newblock {\em Electron. J. Combin.}, 12:Research Paper 46, 12 pp.
  (electronic), 2005.

\bibitem[BKNV05]{bartholdi-k-n-v:bounded}
L.~Bartholdi, Vadim Kaimanovich, V.~Nekrashevych, and B{\'a}lint
Vir{\'a}g.
\newblock Amenability of automata groups.
\newblock (preprint), 2005.

\bibitem[BN]{bond_n:schreier}
I.~Bondarenko and V.~Nekrashevych.
\newblock Growth of {S}chreier graphs of groups generated by bounded automata.
\newblock in preparation.

\bibitem[BN06]{bartholdi_n:rabbit}
Laurent~I. Bartholdi and Volodymyr~V. Nekrashevych.
\newblock Thurston equivalence of topological polynomials.
\newblock {\em Acta Math.}, 197(1):1--51, 2006.

\bibitem[Bos05]{boston:reducing}
Nigel Boston.
\newblock Reducing the {F}ontaine-{M}azur conjecture to group theory.
\newblock In {\em Progress in Galois theory}, volume~12 of {\em Dev. Math.},
  pages 39--50. Springer, New York, 2005.

\bibitem[Bos06]{boston:survey}
Nigel Boston.
\newblock Galois $p$-groups unramified at $p$ - a survey.
\newblock In {\em Primes and Knots}, volume 416 of {\em Contemp. Math.} Amer.
  Math. Soc., 2006.

\bibitem[BS]{bondarenko_s:sushch}
I.~Bondarenko and D.~Savchuk.
\newblock On {S}ushchansky $p$-groups.
\newblock {S}ubmitted.

\bibitem[B{\v{S}}01]{bartholdi-s:wpg}
Laurent Bartholdi and Zoran {\v{S}}uni{\'k}.
\newblock On the word and period growth of some groups of tree automorphisms.
\newblock {\em Comm. Algebra}, 29(11):4923--4964, 2001.

\bibitem[B{\v{S}}06]{bartholdi_s:bsolitar}
Laurent~I. Bartholdi and Zoran {\v{S}}uni{\'k}.
\newblock Some solvable automaton groups.
\newblock In {\em Topological and Asymptotic Aspects of Group Theory}, volume
  394 of {\em Contemp. Math.}, pages 11--29. Amer. Math. Soc., Providence, RI,
  2006.

\bibitem[BV05]{bartholdi-v:basilica}
Laurent Bartholdi and B{\'a}lint Vir{\'a}g.
\newblock Amenability via random walks.
\newblock {\em Duke Math. J.}, 130(1):39--56, 2005.

\bibitem[Con94]{connes:b-ncg}
Alain Connes.
\newblock {\em Noncommutative geometry}.
\newblock Academic Press Inc., San Diego, CA, 1994.

\bibitem[Day57]{day:amenable}
Mahlon~M. Day.
\newblock Amenable semigroups.
\newblock {\em Illinois J. Math.}, 1:509--544, 1957.

\bibitem[Dix50]{dixmier:problem}
Jacques Dixmier.
\newblock Les moyennes invariantes dans les semi-groups et leurs applications.
\newblock {\em Acta Sci. Math. Szeged}, 12(Leopoldo Fejer et Frederico Riesz
  LXX annos natis dedicatus, Pars A):213--227, 1950.

\bibitem[ECH{\etalchar{+}}92]{epstein-al:wp}
David B.~A. Epstein, James~W. Cannon, Derek~F. Holt, Silvio V.~F.
Levy,
  Michael~S. Paterson, and William~P. Thurston.
\newblock {\em Word processing in groups}.
\newblock Jones and Bartlett Publishers, Boston, MA, 1992.

\bibitem[Ers04]{erschler:subexponential}
Anna Erschler.
\newblock Boundary behavior for groups of subexponential growth.
\newblock {\em Ann. of Math. (2)}, 160(3):1183--1210, 2004.

\bibitem[FG91]{fabrykovski-g:growth2}
Jacek Fabrykowski and Narain Gupta.
\newblock On groups with sub-exponential growth functions. {II}.
\newblock {\em J. Indian Math. Soc. (N.S.)}, 56(1-4):217--228, 1991.

\bibitem[FRR95]{ferry-r-r:novikovconj}
Steven~C. Ferry, Andrew Ranicki, and Jonathan Rosenberg.
\newblock A history and survey of the {N}ovikov conjecture.
\newblock In {\em Novikov conjectures, index theorems and rigidity, Vol.\ 1
  (Oberwolfach, 1993)}, volume 226 of {\em London Math. Soc. Lecture Note
  Ser.}, pages 7--66. Cambridge Univ. Press, Cambridge, 1995.

\bibitem[GLS{\.Z}00]{grigorchuk-al:atiyah}
Rostislav~I. Grigorchuk, Peter Linnell, Thomas Schick, and Andrzej
{\.Z}uk.
\newblock On a question of {A}tiyah.
\newblock {\em C. R. Acad. Sci. Paris S\'er. I Math.}, 331(9):663--668, 2000.

\bibitem[Glu61]{glushkov:ata}
V.~M. Glushkov.
\newblock Abstract theory of automata.
\newblock {\em Uspekhi mat. nauk.}, 16(5):3--62, 1961.
\newblock (in Russian).

\bibitem[GNS00]{gns00:automata}
R.~I. Grigorchuk, V.~V. Nekrashevich, and V.~I. Sushchanski{\u\i}.
\newblock Automata, dynamical systems, and groups.
\newblock {\em Tr. Mat. Inst. Steklova}, 231(Din. Sist., Avtom. i Beskon.
  Gruppy):134--214, 2000.

\bibitem[Gre69]{greenleaf:b-means}
Frederick~P. Greenleaf.
\newblock {\em Invariant means on topological groups and their applications}.
\newblock Van Nostrand Mathematical Studies, No. 16. Van Nostrand Reinhold Co.,
  New York, 1969.

\bibitem[Gri80]{grigorchuk:burnside}
R.~I. Grigorchuk.
\newblock On {B}urnside's problem on periodic groups.
\newblock {\em Funktsional. Anal. i Prilozhen.}, 14(1):53--54, 1980.

\bibitem[Gri84]{grigorchuk:gdegree}
R.~I. Grigorchuk.
\newblock Degrees of growth of finitely generated groups and the theory of
  invariant means.
\newblock {\em Izv. Akad. Nauk SSSR Ser. Mat.}, 48(5):939--985, 1984.

\bibitem[Gri00]{grigorchuk:jibg}
R.I. Grigorchuk.
\newblock Just infinite branch groups.
\newblock In Markus P. F. du~Sautoy Dan~Segal and Aner Shalev, editors, {\em
  New horizons in pro-$p$ groups}, pages 121--179. Birkh\"auser Boston, Boston,
  MA, 2000.

\bibitem[GS83a]{gupta_s:pgroups}
N.~Gupta and Said Sidki.
\newblock Some infinite {$p$}-groups.
\newblock {\em Algebra i Logika}, 22(5):584--589, 1983.

\bibitem[GS83b]{gupta-s:burnside}
Narain~D. Gupta and Said~N. Sidki.
\newblock On the {B}urnside problem for periodic groups.
\newblock {\em Math. Z.}, 182(3):385--388, 1983.

\bibitem[G{\v{S}}06]{grigorchuk-s:hanoi-crm}
Rostislav Grigorchuk and Zoran {\v{S}}uni{\'k}.
\newblock Asymptotic aspects of {S}chreier graphs and {H}anoi {T}owers groups.
\newblock {\em C. R. Math. Acad. Sci. Paris}, 342(8):545--550, 2006.

\bibitem[GS{\v{S}}]{grigorch_ss:img}
Rostislav Grigorchuk, Dmytro Savchuk, and Zoran {\v{S}}uni{\'c}.
\newblock The spectral problem, substitutions and iterated monodromy.
\newblock Submitted.

\bibitem[G{\.Z}01]{grigorchuk-z:l2}
Rostislav~I. Grigorchuk and Andrzej {\.Z}uk.
\newblock The lamplighter group as a group generated by a 2-state automaton,
  and its spectrum.
\newblock {\em Geom. Dedicata}, 87(1-3):209--244, 2001.

\bibitem[G{\.Z}02a]{grigorch_z:basilica}
Rostislav~I. Grigorchuk and Andrzej {\.Z}uk.
\newblock On a torsion-free weakly branch group defined by a three state
  automaton.
\newblock {\em Internat. J. Algebra Comput.}, 12(1-2):223--246, 2002.

\bibitem[G{\.Z}02b]{grigorch_z:basilica_sp}
Rostislav~I. Grigorchuk and Andrzej {\.Z}uk.
\newblock Spectral properties of a torsion-free weakly branch group defined by
  a three state automaton.
\newblock In {\em Computational and statistical group theory (Las Vegas,
  NV/Hoboken, NJ, 2001)}, volume 298 of {\em Contemp. Math.}, pages 57--82.
  Amer. Math. Soc., Providence, RI, 2002.

\bibitem[Hin89]{hinz:towers}
Andreas~M. Hinz.
\newblock The {T}ower of {H}anoi.
\newblock {\em Enseign. Math. (2)}, 35(3-4):289--321, 1989.

\bibitem[Ho{\v r}63]{horejs:automata}
Ji{\v r}{\'\i} Ho{\v r}ej{\v s}.
\newblock Transformations defined by finite automata.
\newblock {\em Problemy Kibernet.}, 9:23--26, 1963.

\bibitem[HR04]{holt-r:word}
D.~Holt and C.~R\"over.
\newblock Groups with indexed co-word problem.
\newblock (submitted to Internat. J. Algebra Comput.), 2004.

\bibitem[Kal45]{kalouj:psubgr1}
L.~Kaloujnine.
\newblock Sur le $p$-groupes de sylow du groupe symetriques de degre $p^m$.
\newblock {\em C.R. Acad. Sci. Paris}, 221:222--224, 1945.

\bibitem[Kal48]{kalouj:psubgr2}
L.~Kaloujnine.
\newblock Sur le $p$-groupes de sylow des groupes symetriques finis.
\newblock {\em Ann. Sci. Ecole Norm. Sup.}, 65:239--279, 1948.

\bibitem[Lub94]{lubotzky:expander-book}
Alexander Lubotzky.
\newblock {\em Discrete groups, expanding graphs and invariant measures},
  volume 125 of {\em Progress in Mathematics}.
\newblock Birkh\"auser Verlag, Basel, 1994.

\bibitem[L{\"u}c02]{luck:b-l2}
Wolfgang L{\"u}ck.
\newblock {\em {$L\sp 2$}-invariants: theory and applications to geometry and
  {$K$}-theory}, volume~44 of {\em Ergebnisse der Mathematik und ihrer
  Grenzgebiete. 3. Folge. A Series of Modern Surveys in Mathematics [Results in
  Mathematics and Related Areas. 3rd Series. A Series of Modern Surveys in
  Mathematics]}.
\newblock Springer-Verlag, Berlin, 2002.

\bibitem[Lys85]{lysionok:presentation}
I.~G. Lys{\"e}nok.
\newblock A set of defining relations for the {G}rigorchuk group.
\newblock {\em Mat. Zametki}, 38(4):503--516, 634, 1985.

\bibitem[Mil68]{milnor:problem}
John~W. Milnor.
\newblock Problem 5603.
\newblock {\em Amer. Mat. Monthly}, 75:685--686, 1968.

\bibitem[Nek05]{nekrash:self-similar}
Volodymyr Nekrashevych.
\newblock {\em Self-similar groups}, volume 117 of {\em Mathematical Surveys
  and Monographs}.
\newblock American Mathematical Society, Providence, RI, 2005.

\bibitem[NS04]{nekrash_s:12endomorph}
V.~Nekrashevych and S.~Sidki.
\newblock Automorphisms of the binary tree: state-closed subgroups and dynamics
  of $1/2$-endomorphisms.
\newblock volume 311 of {\em London Math. Soc. Lect. Note Ser.}, pages
  375--404. {Cambridge Univ. Press}, 2004.

\bibitem[Pat88]{paterson:b-amenability}
Alan L.~T. Paterson.
\newblock {\em Amenability}, volume~29 of {\em Mathematical Surveys and
  Monographs}.
\newblock American Mathematical Society, Providence, RI, 1988.

\bibitem[Pil00]{pilgrim:dessins}
Kevin~M. Pilgrim.
\newblock Dessins d'enfants and {H}ubbard trees.
\newblock {\em Ann. Sci. \'Ecole Norm. Sup. (4)}, 33(5):671--693, 2000.

\bibitem[Pis05]{pisier:are}
Gilles Pisier.
\newblock Are unitarizable groups amenable?
\newblock In {\em Infinite groups: geometric, combinatorial and dynamical
  aspects}, volume 248 of {\em Progr. Math.}, pages 323--362. Birkh\"auser,
  Basel, 2005.

\bibitem[Sid87]{sidki:presentation}
Said Sidki.
\newblock On a {$2$}-generated infinite {$3$}-group: the presentation problem.
\newblock {\em J. Algebra}, 110(1):13--23, 1987.

\bibitem[Sid04]{sidki:pol}
Said Sidki.
\newblock Finite automata of polynomial growth do not generate a free group.
\newblock {\em Geom. Dedicata}, 108:193--204, 2004.

\bibitem[Sto94]{stockmeyer:variations}
Paul~K. Stockmeyer.
\newblock Variations on the four-post {T}ower of {H}anoi puzzle.
\newblock {\em Congr. Numer.}, 102:3--12, 1994.

\bibitem[Su{\v{s}}79]{sushchanskii:burnside}
V.~{\=I}. Su{\v{s}}{\v{c}}ans{\cprime}ki{\u\i}.
\newblock Periodic {$p$}-groups of permutations and the unrestricted {B}urnside
  problem.
\newblock {\em Dokl. Akad. Nauk SSSR}, 247(3):557--561, 1979.

\bibitem[vN29]{vneumann:masses}
John von Neumann.
\newblock Zur allgemeinen {Theorie} des {Masses}.
\newblock {\em Fund. Math.}, 13:73--116 and 333, 1929.
\newblock = \emph{Collected works}, vol.\ I, pages 599--643.

\bibitem[VV05]{vorobets:aleshin}
M.~Vorobets and Ya. Vorobets.
\newblock On a free group of transformations defined by an automaton, 2005.
\newblock preprint.

\bibitem[Wag85]{wagon:b-paradox}
Stan Wagon.
\newblock {\em The {B}anach-{T}arski paradox}, volume~24 of {\em Encyclopedia
  of Mathematics and its Applications}.
\newblock Cambridge University Press, Cambridge, 1985.

\bibitem[Wil00]{wilson:jipg}
John~S. Wilson.
\newblock On just infinite abstract and profinite groups.
\newblock In {\em New horizons in pro-$p$ groups}, volume 184 of {\em Progr.
  Math.}, pages 181--203. Birkh\"auser Boston, Boston, MA, 2000.

\bibitem[Wil04]{wilson:nonuniform}
John~S. Wilson.
\newblock On exponential growth and uniformly exponential growth for groups.
\newblock {\em Invent. Math.}, 155(2):287--303, 2004.

\end{thebibliography}

\end{document}